\documentclass[12pt,leqno]{article}
\usepackage{amsmath, amsthm, amsfonts, amssymb, color}
\usepackage{mathrsfs}
\newtheorem{thm}{Theorem}[section]
\newtheorem{cor}[thm]{Corollary}

\theoremstyle{definition}
\newtheorem{defn}{Definition}[section]
\newcommand{\scr}[1]{\mathscr #1}
\definecolor{wco}{rgb}{0.5,0.2,0.3}

\numberwithin{equation}{section} \theoremstyle{remark}
\newtheorem{rem}{Remark}[section]

\newcommand{\ua}{\uparrow}

\title{{\bf Large Deviations for Stochastic Evolution  Equations with
 Small Multiplicative Noise}
\footnote{Supported in part by the DFG through the Internationales
 Graduiertenkolleg
``Stochastics and Real World Models''
and NNSFC(10721091).
} }
\author{{\bf Wei Liu
 \footnote{E-mail: weiliu0402@yahoo.com.cn~~~
Tel:+49-(0)521-1062971, Fax: +49-(0)521-1066455}
}\\
\footnotesize{Fakult\"at F\"ur Mathematik, Universit\"at
Bielefeld, D-33501 Bielefeld, Germany}\\
\footnotesize{School  of Mathematical Sciences, Beijing Normal
 University, Beijing 100875, China}\\
}
\date{}
\begin{document}
\maketitle
\begin{abstract} The Freidlin-Wentzell large deviation principle is
 established for the
 distributions of
  stochastic evolution equations with general monotone drift and
  small
 multiplicative noise.
 As  examples,  the main results are applied to derive the  large deviation principle for
 different types
 of SPDE such as  stochastic reaction-diffusion equations,  stochastic
porous media equations and fast diffusion equations, and the
 stochastic
 $p$-Laplace equation in Hilbert space. The weak convergence approach
 is employed in the proof to establish the Laplace
principle, which is equivalent to the large deviation principle in our
 framework.
\end{abstract} \noindent
 AMS subject Classification:\ 60F10, 60H15.   \\
\noindent
 Keywords: Stochastic evolution equation, large deviation principle, Laplace principle,
variational approach, weak convergence approach, reaction-diffusion
equations, porous media equations, fast diffusion equations,
$p$-Laplace equation.
 \vskip 2cm

\def\R{\mathbb R}  \def\ff{\frac} \def\ss{\sqrt} \def\BB{\mathbb
B}
\def\N{\mathbb N} \def\kk{\kappa} \def\m{{\bf m}}
\def\dd{\delta} \def\DD{\Delta} \def\vv{\varepsilon} \def\rr{\rho}
\def\<{\langle} \def\>{\rangle} \def\GG{\Gamma} \def\gg{\gamma}
  \def\nn{\nabla} \def\pp{\partial} \def\tt{\tilde}
\def\d{\text{\rm{d}}} \def\bb{\beta} \def\aa{\alpha} \def\D{\scr D}
\def\E{\scr E} \def\si{\sigma} \def\ess{\text{\rm{ess}}}
\def\beg{\begin} \def\beq{\begin{equation}}  \def\F{\scr F}
\def\Ric{\text{\rm{Ric}}} \def\Hess{\text{\rm{Hess}}}\def\B{\mathbb B}
\def\e{\text{\rm{e}}} \def\ua{\underline a} \def\OO{\Omega}
 \def\b{\mathbf b}
\def\oo{\omega}     \def\tt{\tilde} \def\Ric{\text{\rm{Ric}}}
\def\cut{\text{\rm{cut}}} \def\P{\mathbb P} \def\ifn{I_n(f^{\bigotimes
 n})}
\def\fff{f(x_1)\dots f(x_n)} \def\ifm{I_m(g^{\bigotimes m})}
 \def\ee{\varepsilon}
\def\pm{\pi_{{\bf m}}}   \def\p{\mathbf{p}}   \def\ml{\mathbf{L}}
 \def\C{\scr C}      \def\aaa{\mathbf{r}}     \def\r{r}
\def\gap{\text{\rm{gap}}} \def\prr{\pi_{{\bf m},\varrho}}
  \def\r{\mathbf r}
\def\Z{\mathbb Z} \def\vrr{\varrho} \def\ll{\lambda}

\def\bt{\begin{theorem}}
\def\et{\end{theorem}}
\def\bl{\begin{lemma}}
\def\el{\end{lemma}}
\def\br{\begin{remark}}
\def\er{\end{remark}}
\def\bx{\begin{Example}}
\def\ex{\end{Example}}
\def\bd{\begin{definition}}
\def\ed{\end{definition}}
\def\bp{\begin{proposition}}
\def\ep{\end{proposition}}
\def\bc{\begin{corollary}}
\def\ec{\end{corollary}}

\newcommand{\ce}{\begin{eqnarray*}}
\newcommand{\de}{\end{eqnarray*}}

\section{Introduction}

There mainly exist three different approaches to analyze
stochastic partial
 differential equations (SPDE) in the literature. The ``martingale measure
 approach" was
 initiated by J.
 Walsh in \cite{Wa}. The ``variational approach" was first
 used by Bensoussan and Temam in \cite{Be,BT} to study SPDE with additive noise, later this approach
 was further developed in the works of Pardoux \cite{Pa}, Krylov and
 Rozovoskii  \cite{KR} for more general case. For the
 ``semigroup (or mild solution)
 approach" we  refer to the classical monograph \cite{DaZa} by Da Prato and Zabcyzk.
In this paper we  use the variational approach to
   treat a large class of  nonlinear SPDE of evolutionary type, which can model
 all kinds of dynamics with stochastic influence in nature or man-made complex systems.
 Stochastic evolution equations
 have been studied intensively  in recent years and we refer to
 \cite{DR1,DRRW,GM,L08a,LW,RRW,R,Wang,Zh}
 for various generalizations and applications.

Concerning the large deviation principle (LDP), there also exist
fruitful
 results within   different frameworks of SPDE.
The general large deviation principle was first formulated by Varadhan \cite{Va1} in 1966. For its
validity
 to stochastic differential equations in  finite dimensional
case we mainly refer to the well known Freidlin-Wentzell LDP
 (\cite{FW}).  The same problem was also treated by
 Varadhan in \cite{Va} and Stroock in \cite{St} by a different
 approach, which followed
 the large deviation theory developed by Azencott \cite{Az},
Donsker-Varadhan \cite{DV} and Varadhan \cite{Va1}.
 In the classical paper \cite{Fr}  Freidlin  studied the large deviations
 for the small noise limit of stochastic reaction-diffusion equations.
 Subsequently, many authors have endeavored to derive the large deviations results
 under less and less restrictive conditions.   We refer the
reader to Da Prato and Zabczyk \cite{DaZa} and Peszat \cite{P} (also
 the references therein) for the extensions to infinite
 dimensional diffusions or stochastic PDE
under global Lipschitz condition on the nonlinear term.
For the case of local Lipschitz conditions we refer to the work of Cerrai and R\"{o}ckner \cite{CR}
where the case of multiplicative and degenerate noise was also investigated. The LDP
for semilinear parabolic equations on a Gelfand triple was studied by Chow in
 \cite{C}. Recently, R\"{o}ckner $et\ al$ established the LDP in
 \cite{RWW} for
 the distributions of the solution to stochastic porous media equations within the variational framework. All these
 papers
 mainly used the classical ideas of discretization approximations and the
 contraction principle, which was first developed by  Freidlin and Wentzell.
 But
 the situation became much involved and complicated in infinite dimensional case
  since each type of  nonlinear SPDE needs different specific
 techniques and estimates.

An alternative approach for LDP has
 been developed by Feng and Krutz in \cite{FK}, which mainly used nonlinear semigroup
 theory and infinite dimensional Hamilton-Jacobi equation. The techniques
rely on the uniqueness theory for the infinite dimensional Hamilton-Jacobi equation
and some exponential tightness estimates.

 In this paper we will study the large deviation principle for
 stochastic evolution
equations with general monotone drift and multiplicative noise,
which
 are more general than the semilinear case studied in \cite{C} and the
 additive noise case in \cite{RWW}. This framework covers all types of SPDE in \cite{R,KR}
such as  stochastic reaction-diffusion equations,  stochastic $p$-Laplace equation,
stochastic porous media equations and fast diffusion equations.
It is quite difficult to follow the classical discretization approach
 in the present case. The reason is many technical difficulties
appear since the coefficients of SPDE in our framework live on a
 Gelfand triple. For example, it is very difficult to obtain some regularity (H\"{o}lder) estimate of
the solution  w.r.t. the time variable, which is essentially required in the classical proof of LDP by discretization approach.

Hence  we would use
 the stochastic control and weak convergence approach in this paper.   This approach
 is mainly based on a variational representation formula for  certain functionals
 of infinite dimensional Brownian Motion, which was established by Budhiraja and Dupuis in \cite{BD}.
The main advantage
 of the weak convergence approach is that one can avoid some exponential probability estimates,
 which might be very difficult to derive for many infinite dimensional
 models. However, in the implement of weak convergence approach, there are still some technical difficulties appearing in the variational framework.
The reason is the coefficients of SEE are nonlinear operators which
are only well-defined via a Gelfand triple (so three
 spaces are involved). Hence we have to properly handle many estimates involving  different spaces instead of just one single space. Some approximation techniques
are also used in the proof.

The weak convergence approach has been used to study the
 large deviations for homeomorphism flows of non-Lipschitz SDEs by Ren and Zhang in \cite{RZ},
 for two-dimensional stochastic Navier-Stokes equations by  Sritharan and Sundar in \cite{SS} and
 reaction-diffusion type SPDEs by Budhiraja $et\ al$ in \cite{BDM}. For more references on
 this approach we may refer to \cite{DE,RZ1,DM08}.

 Let us first recall some standard definitions and results from the large
 deviation theory.
Let $\{X^\varepsilon\}$ be a family of random variables defined on a
 probability space $(\Omega,\mathcal{F},\mathbf{P})$ and taking values
 in some Polish space $E$. Roughly speaking, the large deviation theory
 concerns itself with the exponential decay of the probability measures of
certain kinds of extreme or tail events.
 The rate of such
 exponential decay is expressed by the ``rate function''.

\begin{defn}(Rate function) A function $I: E\to [0,+\infty]$ is called
 a rate function if
$I$ is lower semicontinuous. A rate function $I$ is called a good
rate
 function if  the level
set $\{x\in E: I(x)\le K\}$ is compact for each $K<\infty$.
\end{defn}

\begin{defn}(Large deviation principle) The sequence
 $\{X^\varepsilon\}$ is said to satisfy
 the {\it large deviation principle with  rate function}
 $I$ if  for each
Borel subset $A$ of $E$
$$
 -\inf_{x\in A^o} I(x)\le \liminf_{\vv\to 0}
   \vv^2 \log \mathbf{P}(X^{\vv}\in A)\le \limsup_{\vv\to 0}
   \vv^2 \log \mathbf{P}(X^{\vv}\in A)\le
  -\inf_{x\in \bar A} I(x),
$$
where $A^o$ and $\bar A$ are respectively the interior  and the
closure of $A$ in $E$.
\end{defn}

If one is interested in obtaining the exponential estimates on
general
 functions instead of the indicator
functions of Borel sets in $E$, then one can study the
 following Laplace principle (LP).

\begin{defn}(Laplace principle) The sequence $\{X^\varepsilon\}$ is
 said to satisfy
 the  {\it Laplace principle  with  rate function} $I$ if
  for each
bounded continuous real-valued function $h$ defined on $E$
$$\lim_{\vv\to 0}\vv^2 \log \mathbf{E}\left\lbrace
 \exp\left[-\frac{1}{\varepsilon^2} h(X^{\vv})\right]\right\rbrace
= -\inf_{x\in E}\left\{h(x)+I(x)\right\}.$$
\end{defn}

The starting point for the weak convergence approach is the
equivalence
 between LDP and LP if
$E$ is a Polish space and the rate function is good. This result was first
formulated in \cite{Pu} and it is essentially a consequence of Varadhan's lemma
 \cite{Va1} and
 Bryc's converse theorem \cite{Br}. We refer to \cite{DE,DZ} for an elementary proof of it.

Let $\{W_t\}_{t\geq0}$ be a cylindrical Wiener process on a
separable
 Hilbert space $U$
w.r.t a complete filtered probability space
$(\Omega,\mathcal{F},\mathcal{F}_t,\mathbf{P})$ (i.e. the path of
 $W$ take
values
 in  $C([0,T];U_1)$, where $U_1$ is another Hilbert space such that the
 embedding $U\subset U_1$ is
 Hilbert-Schmidt). Suppose $g^\varepsilon: C([0,T]; U_1)\rightarrow E$
 is a measurable
map and
 $X^\varepsilon= g^\varepsilon(W_{\cdot})$.
Let
$$\mathcal{A}=\left\lbrace v: v\  \text{is  $U$-valued
 $\mathcal{F}_t$-predictable process and}\
  \int_0^T\|v_s(\omega)\|^2_U\d s<\infty \  a.s.\right\rbrace, $$
 $$S_N=\left\lbrace \phi\in L^2([0,T], U):
\int_0^T\|\phi_s\|^2_{U}\d
 s\leq N
 \right\rbrace.$$ The set $S_N$ endowed with the weak topology is a
Polish space (we will always refer to the weak topology on $S_N$ in this
 paper if we don't state it explicitly).  Define
 $$\mathcal{A}_N=\left\{v\in\mathcal{A}: v(\omega)\in S_N\
 \mathbf{P}-a.s.\right\}.$$
 Now
we  formulate the following sufficient condition for the Laplace
 principle (equivalently, large deviation principle) of $X^\vv$ as
 $\vv\rightarrow0$.\\

$\textbf{(A)}$  There exists a measurable map
$g^0:
 C([0,T]; U_1)\rightarrow E$
  such that the following
two conditions hold:

(i) Let $\{v^\varepsilon: \varepsilon>0\}\subset \mathcal{A}_N$ for
 some $N<\infty$. If $v^\varepsilon$ converge to $v$ in distribution
  as
 $S_N$-valued random elements, then
 $$g^\varepsilon\left(W_\cdot+\frac{1}{\varepsilon}
\int_0^\cdot v^\varepsilon_s\d s \right)\rightarrow
g^0\left(\int_0^\cdot
 v_s\d s \right)$$
 in distribution as $\vv\rightarrow 0$.

(ii) For each $N<\infty$, the set
$$K_N=\left\{g^0\left(\int_0^\cdot \phi_s\d
s\right):
 \phi\in S_N\right\}$$
  is a compact subset of $E$.

\beg{lem}\label{L1}\cite[Theorem 4.4]{BD}  If $X^\varepsilon=g^\varepsilon(W)$ and
 the assumption $\bf{(A)}$ holds, then the family
  $\{X^\varepsilon\}$ satisfies the Laplace principle (hence large
 deviation
 principle) on $E$ with the good
rate function $I$ given by
\begin{equation}\label{rate formula}
I(f)=\inf_{\left\{\phi\in L^2([0,T]; U):\  f=g^0(\int_0^\cdot \phi_s\d
 s)\right\}}\left\lbrace\frac{1}{2}
\int_0^T\|\phi(s)\|_U^2\d s \right\rbrace.
\end{equation}
\end{lem}

We will verify the sufficient condition $\textbf{(A)}$ for general SPDE within the variational framework.
Besides the classical monotone conditions assumed for the well-posedness of SPDE, we need to
require one additional assumption (see $(A4)$ below) on the noise coefficient for the LDP.
In fact, the weak convergence approach are used here to avoid the  time discretization for SPDE (the most technical
and difficult step in the classical proof of LDP)
since the regularity estimate of the solution w.r.t. the time variable is unavailable in the variational framework.
But unlike the semilinear case (e.g.\cite{BDM}),  we have to  use  It\^{o}'s formula for the square norm
of the solution in the estimate. Then the weak convergence of control $v^\varepsilon$ to $v$ (see (i) of $\textbf{(A)}$) cause some technical difficulty in the proof of convergence of corresponding solutions under the variational framework. Hence  we need to have some restriction on the noise (see $(A5)$) such that the weak convergence procedure can be verified. Later some standard approximation techniques are used to relax this assumption.

\section{Main framework and result}
 Let
$$V\subset H\equiv H^*\subset V^*$$
 be a Gelfand triple, i.e. $V$
is a reflexive and separable Banach space and $V^*$ is its dual space, $(H,\<\cdot,\cdot\>_H)$ is a separable
Hilbert space and identified with its dual space by Riesz
isomorphism, $V$ is
 continuously and densely embedded in $H$. The dualization
between $V^*$ and $V$  is denoted by $_{V^*}\<\cdot,\cdot\>_V$ and it is
 obvious that
$${ }_{V^*}\<u, v\>_V=\<u, v\>_H, \  u\in H ,v\in V.$$
Let $\{W_t\}_{t\geq0}$ be a cylindrical Wiener process on a separable
 Hilbert space $U$
w.r.t a complete filtered probability space
$(\Omega,\mathcal{F},\mathcal{F}_t,\mathbf{P})$. $\left(L_2(U;H)
 \|\cdot\|_2\right)$ denote
the space of all  Hilbert-Schmidt operators  from $U$ to $H$. We use
$L(X,Y)$ to denote the space of all bounded linear operators from space $X$ to
$Y$.

 Consider
the following stochastic evolution equation \beq\label{**}
  \d X_t = A(t,X_t)\d t+B(t,X_t)\d W_t,
\end{equation}
where $A: [0,T]\times V\to V^*$ and $B: [0,T]\times V \to L_2(U;H)$
are
 measurable. For the large deviation
principle we need to assume the following conditions, which are slightly
stronger than those assumed in \cite{KR} for the existence and uniqueness of
strong solution to (\ref{**}).

For a fixed $\alpha>1$, there exist constants $\delta>0$ and $K$
such that the
 following
 conditions hold for all $v,v_1,v_2\in V$ and $t\in [0,T]$.
\begin{enumerate}
\item [$(A1)$] (Hemicontinuity) The map
$ s\mapsto  { }_{V^*}\<A(t,v_1+s v_2),v\>_V$
     is continuous on $\mathbb{R}$.
\item [$(A2)$] (Strong monotonicity)
$$2{ }_{V^*}\<A(t,v_1)-A(t,v_2), v_1-v_2\>_V+
 \|B(t,v_1)-B(t,v_2)\|_{2}^2
    \leq -\delta\|v_1-v_2\|_V^\alpha+ K\|v_1-v_2\|_H^2.$$
\item[$(A3)$] (Boundedness) $\sup_{t\in[0,T]}\|B(t,0)\|_2<\infty$ and
 $$  \|A(t,v)\|_{V^*}+ \|B(t,v)\|_{L(U,V^*)}\leq
 K(1+\|v\|_V^{\alpha-1}).$$
\item[$(A4)$]
Suppose
 there exist a sequence of subspaces $\{H_n\}$ such that
 $$H_n\subseteq H_{n+1}, \  H_n\hookrightarrow V \text{compact and}\  \bigcup_{n=1}^\infty H_n\subseteq H \
 \text{dense},$$
 and for any $M>0$
 \begin{equation}\label{approximated property}
 \sup_{(t,v)\in
 [0,T]\times S_M}\|P_nB(t,v)-B(t,v)\|_2\rightarrow 0\ (n\rightarrow\infty),
 \end{equation}
 where $P_n: H\rightarrow H_n$ is the projection operator and
 $S_M=\{v\in V: \|v\|_H\le M \}$.
\end{enumerate}

 \beg{rem}\label{rem1}(i) By $(A2)$ and $(A3)$ we can easily obtain the
coercivity and boundedness of $A$ and $B$:
$$2_{V^*}\<A(t,v),
 v\>_V+\|B(t,v)\|_{2}^2+\frac{\delta}{2}\|v\|_V^\alpha
    \le C(1+ \|v\|_H^2),$$
$$\|B(t,v)\|_2^2\le C(1+\|v\|_H^2+\|v\|_V^\alpha).$$
Hence the boundedness of $B$ in $(A3)$ automatically holds if $\alpha\ge 2$. If $1<\alpha<2$, the additional
assumption on $B$ in $(A3)$ is assumed for the well-posedness of the skeleton equation (see (\ref{1.4})).

 (ii)
 Since for all $(t,v)\in [0,T]\times V$ we have
$$\|P_nB(t,v)-B(t,v)\|_2\rightarrow 0\ (n\rightarrow\infty). $$
Hence a simple sufficient condition
for (\ref{approximated property}) holds is to assume that
$$\left\{B(t,v): (t,v)\in
 [0,T]\times S_M  \right\}$$
 is a relatively compact set in $L_2(U;H)$. For example, we can take
 $$B(t,v)=\sum_{i=1}^N b_i(v)B_i(t), $$
 where $b_i(\cdot): V\rightarrow \mathbb{R}$ are Lipschitz functions and
 $B_i(\cdot):[0,T]\rightarrow L_2(U;H)$ are continuous.

Another simple example is $B(t,v)=QB_0(t,v)$ where $Q\in L_2(H; H)$ and
$$B_0: [0,T]\times V\rightarrow L(U;H),\ \
\sup_{(t,v)\in [0,T]\times S_M}\|B_0(t,v)\|_{L(U;H)}< \infty, \ \forall M>0.$$

 (iii) If there exists a Hilbert space $H_0$ such that the
embedding $H_0\subseteq H$ is compact, $\{e_i\}\subseteq H_0\cap V$
is an ONB in $H_0$ and also orthogonal  in $H$. Suppose for all
$M>0$
$$\sup_{(t,v)\in [0,T]\times S_M}\|B(t,v)\|_{L_2(U;H_0)}< \infty.$$
Then $(\ref{approximated property})$ holds. Because
$B(t,v)=\sum_{i,j=1}^\infty b_{i,j}(t,v)u_i\otimes e_j$, by
assumptions we know $\|e_j\|_H^2\rightarrow 0$ and
$$\sup_{(t,v)\in [0,T]\times S_M}\sum_{i,j=1}^\infty b_{i,j}^2(t,v)<\infty.$$
then
$$\|P_nB(t,v)-B(t,v)\|_2^2=\sum_{i=1}^\infty \sum_{j=n+1}^\infty
b_{i,j}^2(t,v)\|e_j\|_H^2.$$
Hence $(\ref{approximated property})$
follows from the dominated convergence theorem.  \qed


\end{rem}

If $(A1)-(A3)$ hold, according to \cite[Theorem II2.1]{KR}
  for any
$X_0\in L^2(\Omega\to H; \mathcal{F}_0;\mathbf{P})$
    (\ref{**})
    has an unique solution $\{X_t\}_{t\in [0,T]}$ which is an adapted
 continuous process
    on $H$ such that $\mathbf{E}\int_0^T\left(
 \|X_t\|_V^{\alpha}+\|X_t\|_H^2\right) \d t<\infty$ and
$$\<X_t, v\>_H= \<X_0,v\>_H +\int_0^t { }_{V^*}\<A(s,X_s),
    v\>_V\d s + \int_0^t\<B(s,X_s)\d W_s, v\>_H, \ \mathbf{P}-a.s.$$
holds for all $v\in V$ and $ t\in
    [0,T]$.  Moreover, we have
$\mathbf{E}\sup_{t\in[0,T]}\|X_t\|_H^2 <\infty$ and the  crucial
It\^{o} formula
$$\|X_t\|_H^2= \|X_0\|_H^2 +\int_0^t\left(2{ }_{V^*}\<A(s,X_s),
    X_s\>_V+\|B(s,X_s)\|_2^2 \right)  \d s + 2\int_0^t\<X_s,
  B(s,X_s)\d W_s\>_H.$$
Let us  consider the general stochastic evolution equation  with
small noise: \beq\label{1.3} \d X_t^{\vv} = A(t,X_t^\vv )\d t+\vv
B(t,X_t^\vv)\d
 W_t,\ \ \
\vv>0,\  X_0^\vv =x\in H.\end{equation}
Hence the unique strong solution $\{X^\vv\}$ of (\ref{1.3}) takes values in
 $C([0,T];
 H)\cap L^\alpha([0,T]; V)$. It's well-known that $\left(C([0,T]; H)\cap L^\alpha([0,T]; V), \rho\right)$ is a Polish
 space with the following metric
\beq\label{metric}\rho(f,g):=\sup_{t\in[0,T]}\|f_t-g_t\|_H
+\left(\int_0^T\|f_t-g_t\|_V^\alpha\d t \right)^{\frac{1}{\alpha}}.
\end{equation}
It follows (from infinite dimensional version
 of Yamada-Watanabe theorem in  \cite{RSZ}) that there exists
 a
 Borel-measurable function
$$g^\varepsilon: C([0,T]; U_1)\rightarrow
 C([0,T]; H)\cap L^\alpha([0,T]; V)$$
such that
 $X^\varepsilon=g^\varepsilon (W)\ a.s.$.  To state our main
result, let us  introduce the skeleton equation associated to
(\ref{1.3}): \beq\label{1.4} \ff{\d z_t^\phi }{\d t} = A(t,z_t^\phi)
+
 B(t,z_t^\phi)\phi_t,\
\ \ z_0^\phi= x,\ \phi\in L^2([0,T];U).
\end{equation}
  An
element $z^\phi \in C([0,T]; H) \cap L^{\alpha}([0,T]; V)$ is called
a solution to (\ref{1.4}) if for any $v\in V$
 \beq\label{1.5}
\<z_t^\phi , v\>_H = \<x,v\>_H + \int_0^t \ _{V^*}\< A(s,z_s^\phi
)+B(s,z_s^\phi)\phi_s, v\>_V \d s,\ \ \ t\in
 [0,T].\end{equation}
We will prove (see Lemma \ref{L2.1}) that $(A1)-(A3)$ also imply the
 existence and uniqueness of
the solution to (\ref{1.4}) for any $\phi\in L^2([0,T];U)$.

Define $g^0: C([0,T]; U_1)\rightarrow C([0,T]; H)\cap
L^\alpha([0,T];
 V)$ by
$$g^0(h):= \beg{cases} z^\phi, &\text{if}\ \  h=\int_0^\cdot \phi_s\d s\ \
 \text{for\ some}\  \phi\in L^2([0,T];U);\\
0, &\text{otherwise.}\end{cases}$$
Then it's obvious that the rate function in (\ref{rate formula}) can
be written as
 \beq\label{rate} I(z)=\inf
\left\{\frac{1}{2}\int_0^T\|\phi_s\|_{U}^2\d s:\
 z=z^{\phi}
,\ \phi\in L^2([0,T], U)\right\},
\end{equation}
where  $z\in C([0,T];H)\cap
 L^\alpha([0,T]; V)$.

Now we  formulate the main result which is a
 Freidlin-Wentzell type estimate.

\beg{thm}\label{T1.1} Assume $(A1)-(A4)$ hold. For each $\vv>0$, let
$X^\vv=\{X^\vv_t\}_{t\in [0,T]}$ be the solution to $(\ref{1.3})$.
Then as $\vv\to 0$, $\{X^\vv\}$ satisfies the $LDP$ on $C([0,T];
H)\cap L^\alpha([0,T]; V)$
 with the good rate function $I$ which is given by $(\ref{rate})$.
\end{thm}

\beg{rem} (i) According to \cite[Theorem 5]{BDM}, we can also prove
uniform
 Laplace principle by using
the same arguments but with more cumbersome notation.

(ii) This theorem can not be applied to stochastic fast-diffusion
 equations in \cite{LW,RRW} since
$(A2)$ fails to satisfy. However,
 if we  replace $(A2)$
by the classical monotone and coercive conditions in \cite{KR}
\begin{enumerate}
\item [$(A2^\prime)$] $$\aligned
2{ }_{V^*}\<A(t,v_1)-A(t,v_2), v_1-v_2\>_V +
 \|B(t,v_1)-B(t,v_2)\|_{2}^2
    & \leq  K\|v_1-v_2\|_H^2,\\
2{ }_{V^*}\<A(t,v), v\>_V+
 \|B(t,v)\|_{2}^2+\delta\|v\|_V^\alpha
    & \leq  K(1+\|v\|_H^2).
\endaligned$$
\end{enumerate}
Then the LDP
can be established
on $C([0,T]; H)$  by the
similar and simpler argument.
\end{rem}

\beg{thm}\label{T1.2} Assume $(A1),(A2^\prime),(A3)-(A4)$ hold.
Then as $\vv\to 0$, the solution $\{X^\vv\}$ of  $(\ref{1.3})$ satisfies the $LDP$ on $C([0,T]; H)$
 with the good rate function $I$ which is given by $(\ref{rate})$.
\end{thm}
\begin{rem} Note that $(A2)$ mainly used to prove the additional
convergence in $L^\alpha([0,T]; V)$.
 Hence, if we only concern the LDP on $C([0,T];H)$, then we can prove
the  Theorem \ref{T1.2} under the  weaker
assumptions above.
Since the proof is  only a small modification (only consider the convergence in
$C([0,T];H)$) of the argument for Theorem
 \ref{T1.1}, we omit the details here.
\end{rem}

The organization of the paper is as follows. In section 3, under the
 additional assumption $(A5)$ on $B$ we  prove
Theorem \ref{T1.1} by using the weak convergence approach.
 Section 4 is devoted to relax
 the assumption $(A5)$  by some standard approximation techniques.
 In section 5 we
 apply the main results to different class of SPDEs in Hilbert space
 as applications.

\section{Proof of Theorem \ref{T1.1} under additional assumption}

In order to verify the sufficient conditions $\textbf{(A)}$, we
 need to first consider the finite dimensional noise, i.e. we approximate the diffusion coefficient $B$
by $P_nB$. But for the simplicity of the notation, we  formulate the following additional assumption on $B$:
\begin{enumerate}
   \item [$(A5)$]  $B:[0,T]\times V\rightarrow L(U;V_0)$ satisfies
$$\|B(t,v)\|_{L(U;V_0)}^2\le C(1+\|v\|_V^\alpha+\|v\|_H^2), $$
\noindent where  $V_0\subseteq
 V$ is compact embedding and $C>0$ is a constant.
\end{enumerate}

For the reader's convenience,  we recall two well-known inequalities
which
 used quite often in the proof. Throughout the paper, the generic
 constants may be different from line to line. If it is essential, we
 will
 write the dependence of the constant on parameters explicitly.

\noindent\textbf{Young's inequality:} Given $p,q>1$ satisfy
 $\frac{1}{p}+\frac{1}{q}=1$, then for any positive number $\sigma, a, b$ we have
$$ab\le \sigma\frac{ a^p}{p}+\sigma^{-\frac{q}{p}}\frac{b^q}{q}.$$

\noindent\textbf{Gronwall's lemma:} Let $F,\Phi,\Psi:
 [0,T]\rightarrow \mathbb{R}^+$ be Lebesgue measurable. Suppose $\Psi$
 is locally
 integrable and $\int_0^T\Psi(s)F(s)\d s<\infty$. If
 \beq\beg{split}
F(t)& \le \Phi(t)+\int_0^t\Psi(s)F(s)\d s, \ t\in[0,T]\ \
  \text{or}\\
\frac{\d}{\d t}F(t)& \le \frac{\d}{\d t}\Phi(t)+\Psi(t)F(t), \
 t\in[0,T),  \ F(0)\le\Phi(0).
\end{split}\end{equation}
Then \beq
 F(t) \le \Phi(t)+\int_0^t\exp\left[ \int_s^t\Psi(u) \d
 u\right]\Psi(s)\Phi(s) \d s, \ t\in[0,T].
\end{equation}

\beg{lem} \label{L2.1} Assume $(A1)-(A3)$ hold. Let
$$\|z\|:=\sup_{t\in [0,T]}\|z_t\|_H^2+\delta\int_0^T\|z_t\|_V^\alpha \d t$$
for $z\in C([0,T];H)\cap L^\alpha([0,T];V)$.
 For all $x\in H$ and $\phi\in L^2([0,T];U)$ there
 exists a unique solution $z^\phi$
 to $(\ref{1.4})$ and
 \beq\label{stable}\beg{split}   \|z^\phi -z^\psi\|
  \le\exp\left\lbrace
\int_0^T\left(K+\|\phi_t\|_U^2+\|B(t,z_t^\psi)\|_2^2 \right)\d t
 \right\rbrace\int_0^T \|\phi_t -\psi_t\|_U^2\d t
 \end{split}
\end{equation}
 hold for some constant $K$ and all $ \phi,\psi\in
 L^2([0,T];U)$.
\end{lem}

\beg{proof} To verify the existence of the solution, we  make use of
\cite[Theorem II.2.1]{KR}. First
we assume $\phi\in
 L^\infty([0,T];U)$ and
$$\tilde{A}(s,v):= A(s,v) + B(s,v)\phi_s.$$
Then, due to $(A1)-(A3)$,  it's easy to verify that $\tilde{A}$
 satisfies Assumptions
$A_i) (i=1,..,5)$ on page 1252 of \cite{KR}.

(i) Hemicontinuity of $\tilde{A}$ follows from $(A1)$ and $(A2)$.

(ii) Monotonicity and coercivity of $\tilde{A}$ follows from $(A2)$ and $(A3)$.

(iii) Boundedness of $\tilde{A}$ follows from $(A3)$.

 Therefore, by \cite{KR}(or \cite[Theorem
30.A]{Z})we know (\ref{1.4}) has an unique solution.

For general $\phi\in L^2([0,T];U)$, we can find a sequence of
$\phi^n\in L^\infty([0,T];U)$ such that
$$\phi_n\rightarrow \phi \ \ \text{strongly in}\ \
 L^2([0,T];U).$$
 Let $z^n$ be the unique solution to
(\ref{1.4}) for $\phi^n$, we will show $\{z^n\}$
 is a Cauchy sequence in $C([0,T];H)\cap L^{\alpha}([0,T];V)$.
By using
$(A2)$ we
 have
\beq\label{cauchy}\beg{split}
 \ff{\d}{\d t} \|z_t^n -z_t^m\|_H^2 =&
2 { }_{V^*} \< A(t,z_t^n )
-A(t,z_t^m), z_t^n-z_t^m\>_V\\
& +2 \< B(t,z_t^n)\phi^n_t
-B(t,z_t^m)\phi^m_t, z_t^n-z_t^m\>_H\\
 \le& 2 { }_{V^*} \< A(t,z_t^n )
-A(t,z_t^m), z_t^n-z_t^m\>_V+ \|B(t,z_t^n)
-B(t,z_t^m)\|_{2}^2\\
& +\|\phi^n_t\|_U^2\|z_t^n -z_t^m\|_H^2
 +2\<z_t^n-z_t^m, B(t,z_t^m)\phi^n_t -B(t,z_t^m)\phi^m_t\>_H\\
\le & -\delta \|z_t^n -z_t^m\|_V^\alpha +
 (K+\|\phi^n_t\|_U^2)\|z_t^n-z_t^m\|_H^2\\
& +2\|B^*(t,z_t^m)\left(z_t^n-z_t^m\right)\|_U\|\phi^n_t
 -\phi^m_t\|_U\\
\le & -\delta \|z_t^n -z_t^m\|_V^\alpha + \|\phi^n_t -\phi^m_t\|_U^2
\\
& + \left(K+\|\phi^n_t\|_U^2+\|B(t,z_t^m)\|_2^2\right)\|z_t^n
-z_t^m\|_H^2 .
\end{split}\end{equation}
where $B^*$ denote the adjoint operator of $B$ and we also use the
fact
$$  \|B^*\|_{L(H;U)}=\|B\|_{L(U;H)}\le\|B\|_2 .$$
Then by the Gronwall lemma we have \beq\label{estimate}\beg{split}
\|z^n -z^m\|
 & \le\exp\left\lbrace
\int_0^T\left(K+\|\phi^n_t\|_U^2+\|B(t,z_t^m)\|_2^2 \right)\d t
 \right\rbrace\int_0^T \|\phi^n_t -\phi^m_t\|_U^2\d t.
 \end{split}
\end{equation}
By the similar argument we have \beq\label{bound}\beg{split}
 \ff{\d}{\d t} \|z_t^n\|_H^2 =&
2{ }_{V^*} \< A(t,z_t^n )
, z_t^n\>_V +2 \< B(t,z_t^n)\phi^n_t, z_t^n\>_H\\
\le & -\frac{\delta}{2} \|z_t^n\|_V^\alpha + C(1+ \|z_t^n\|_H^2)+
 \|\phi^n_t\|_U^2\|z_t^n\|_H^2.
\end{split}\end{equation}
Then by the Gronwall lemma and boundedness of $\phi^n$  in
 $L^2([0,T]; U)$
\beq\label{bound of z^n}
 \|z^n \| \le C\exp\left\lbrace
\int_0^T\left(C+\|\phi^n_t\|_U^2\right)\d t \right\rbrace \left(
 \|x\|_H^2+T\right)\le \textbf{Constant}< \infty.
\end{equation}
Hence we have \beq\label{bound of B} \int_0^T\|B(t,z_t^m)\|_2^2\d t
\le
 C\int_0^T\left(1+\|z_t^m\|_H^2+\|z_t^m\|_V^\alpha\right)\d t\le
 \textbf{Constant}< \infty.
\end{equation}
Combining (\ref{estimate}),(\ref{bound of B}) and
 $\phi^n\rightarrow \phi$, we can conclude that
$\{z^n\}$ is a Cauchy sequence in $C([0,T];H)\cap
L^{\alpha}([0,T];V)$, and we denote
 the limit by $z^\phi$.

Then by repeating the standard monotonicity argument(e.g.\cite[Theorem 30.A]{Z}) one
 can show that $z^\phi$ is the solution of (\ref{1.4})
 corresponding to $\phi$.

And (\ref{stable}) can be derived from
(\ref{estimate}). Hence the proof
 is complete.
\end{proof}

The following result shows that $I$ defined by (\ref{rate}) is a good
 rate function.
\beg{lem}\label{L2.2} Assume $(A1)-(A3)$ hold. For every $N<\infty$,
the set
$$K_N=\left\{g^0\left(\int_0^\cdot \phi_s\d s\right):
 \phi\in S_N\right\}$$
  is a compact subset in $C([0,T]; H)\cap
 L^\alpha([0,T];V)$.
\end{lem}
\begin{proof}\textbf{Step 1}: we first assume $B$ also satisfy $(A5)$.
 By definition we know
$$K_N=\left\{z^\phi: \phi\in L^2([0,T];U),\  \int_0^T\|\phi_s\|_{U}^2\d s\le
 N \right\}.$$
For any sequence ${\phi^n}\subset S_N$, we may assume
 $\phi^n\rightarrow\phi$ weakly in $L^2([0,T];U)$ since  $S_N$
 is weakly
 compact. Denote $z^n$ and $z$ are the solutions of (\ref{1.4})
 corresponding to $\phi^n$ and $\phi$ respectively. Now it's sufficient  to
 show
 $z^n\rightarrow z$ strongly in $C([0,T]; H)\cap L^\alpha([0,T];V)$.

 From
 (\ref{cauchy}) we have
\begin{equation}
\begin{split}\label{difference}
& \|z^n_t -z_t\|_H^2+\delta\int_0^t\|z^n_s -z_s\|_V^\alpha\d s\\
&  \le \int_0^t(K+\|\phi^n_s\|_U^2)\|z_s^n -z_s\|_H^2\d s +2\int_0^t
\<z_s^n-z_s, B(s,z_s)(\phi^n_s -\phi_s)\>_H\d s.
\end{split}
\end{equation}
Define
$$h_t^n=\int_0^t B(s,z_s)(\phi^n_s
-\phi_s)\d s . $$ By $(A5)$ and (\ref{bound of B}) we know $h^n\in
 C([0,T];V_0)$ and
\begin{equation}\begin{split}\label{norm estimate}
\sup_{t\in[0,T]}\|h^n_t\|_{V_0}& \le \int_0^T\| B(s,z_s)(\phi^n_s
 -\phi_s)\|_{V_0}\d s\\
& \le \left( \int_0^T\| B(s,z_s)\|_{L(U,V_0)}^2\d s\right)
 ^{1/2}\left(\int_0^T\| \phi^n_s -\phi_s\|_{U}^2\d s \right)^{1/2} \\
& \le \textbf{Constant}< \infty.
\end{split}
\end{equation}
Since  the embedding $V_0\subseteq V$ is compact and
 $\phi^n\rightarrow\phi$ weakly in $L^2([0,T];U)$, it's easy to show
 that  $h^n\rightarrow
 0$ in $C([0,T]; V)$ by using the Arz\`ela-Ascoli theorem(also see
 e.g.\cite[Lemma 3.2]{BD})
 (more precisely, this convergence may only hold for a subsequence,
  but it's enough for our purpose since we may denote the convergent
 subsequence still by $h^n$).
 In particular,  $h^n\rightarrow 0$ in $C([0,T]; H)\cap
 L^\alpha([0,T]; V)$.

 Moreover the derivative (w.r.t. time variable) is given by
 $$(h_s^n)^\prime= B(s,z_s)(\phi^n_s -\phi_s).$$
As in the Lemma \ref{L2.1}, we may assume $\phi^n, \phi\in L^\infty([0,T]; U)$ first.
Then by $(A3)$
\begin{equation}\begin{split}
\int_0^T\|(h^n_s)^\prime\|^{\frac{\alpha}{\alpha-1}}_{V^*}\d s& \le
 \int_0^T\| B(s,z_s)(\phi^n_s
 -\phi_s)\|^{\frac{\alpha}{\alpha-1}}_{V^*}\d
 s\\
& \le
C\int_0^T\left(1+\|z_s\|_{V}^{\alpha} \right)\d s \\
& \le \textbf{Constant}< \infty.
\end{split}
\end{equation}
Hence  $(h_\cdot^n)^\prime$ is an element in
 $L^{\frac{\alpha}{\alpha-1}}([0,T]; V^*)$.

By \cite[Proposition 23.23]{Z} we have the following integration by
parts
 formula
$$\<z_t^n-z_t, h^n_t\>_H=\int_0^t{ }_{V^*}\<(z_s^n-z_s)^\prime,
 h^n_s\>_V\d s+\int_0^t{ }_{V^*}\< (h_s^n)^\prime, z_s^n-z_s\>_V\d s.$$
 Hence one has
\begin{equation}\begin{split}\label{e1}
& \int_0^t \<z_s^n-z_s, B(s,z_s)(\phi^n_s -\phi_s)\>_H\d s\\
=& \<z_t^n-z_t, h^n_t\>_H- \int_0^t{ }_{V^*}\<(z_s^n-z_s)^\prime,
 h^n_s\>_V\d s\\
=&  \<z_t^n-z_t, h^n_t\>_H-\int_0^t{ }_{V^*}\<A(s,z_s^n)-A(s,z_s),
 h^n_s\>_V\d s\\
 & -\int_0^t \<B(s,z_s^n)\phi_s^n-B(s,z_s)\phi_s,
 h^n_s\>_H\d s\\
=:& I_1+I_2+I_3
\end{split}
\end{equation}
By using the H\"{o}lder inequality, $(A3)$ and (\ref{bound of z^n})
we
 have
\beq\label{e2}
\begin{split}
I_1& \le \|z_t^n-z_t\|_H\cdot\|h^n_t\|_H\le
 \frac{1}{4}\|z_t^n-z_t\|_H^2+\|h^n_t\|^2_H .\\
I_2 & \le \int_0^t \|A(s,z_s^n)-A(s,z_s)\|_{V^*}\|h^n_s\|_V\d s\\
& \le \left( \int_0^t
 \|A(s,z_s^n)-A(s,z_s)\|_{V^*}^{\frac{\alpha}{\alpha-1}}\d
 s\right)^{\frac{\alpha-1}{\alpha}}
\left( \int_0^t \|h^n_s\|_V^\alpha\d s\right)^{\frac{1}{\alpha}} \\
& \le \left(\int_0^t C\left(
 1+\|z_s\|_V^{\alpha}+\|z_s^n\|_V^{\alpha}\right)
\d s\right)^{\frac{\alpha-1}{\alpha}}
\left( \int_0^t \|h^n_s\|_V^\alpha\d s\right) ^{\frac{1}{\alpha}}\\
& \le C \left( \int_0^t \|h^n_s\|_V^\alpha\d s\right)
 ^{\frac{1}{\alpha}}.\\
I_3 & \le \int_0^t \|B(s,z_s^n)\phi_s^n-B(s,z_s)\phi_s\|_H\cdot\|
 h^n_s\|_H\d s\\
& \le \sup_{s\in[0,t]}\| h^n_s\|_H
\int_0^t \|B(s,z_s^n)\phi_s^n-B(s,z_s)\phi_s\|_H\d s\\
& \le \sup_{s\in[0,t]}\| h^n_s\|_H \left\lbrace N^{1/2}\left(
\int_0^t \|B(s,z_s^n)\|_2^2\d s\right)^{1/2} + N^{1/2}\left(
\int_0^t \|B(s,z_s)\|_2^2\d s\right)^{1/2}
\right\rbrace \\
& \le C \sup_{s\in[0,t]}\| h^n_s\|_H.
\end{split}
\end{equation}
where $C$ is a constant which come from the following estimate
$$\int_0^t\|B(s,z_s^n)\|^2_{2}\d s\le
 C\int_0^t\left(1+\|z_s^n\|^2_{H}+ \|z_s^n\|^\alpha_{V}\right) \d s
 \le \textbf{Constant}<\infty $$
Combining (\ref{difference}) and (\ref{e1})-(\ref{e2})  we have
\begin{equation}\begin{split}
& \|z^n_t -z_t\|_H^2+\delta \int_0^t \|z^n_t -z_t\|_V^\alpha\d t\\
& \le C\int_0^t(1+\|\phi^n_s\|_U^2)\|z_s^n -z_s\|_H^2\d s
 +C\left(\sup_{s\in[0,t]}\| h^n_s\|_H +\sup_{s\in[0,t]}\|
h^n_s\|_H^2 + \left( \int_0^t \|h^n_s\|_V^\alpha\d
s\right)^{\frac{1}{\alpha}}
 \right)
\end{split}\end{equation}
Then by the Gronwall lemma and
 $L^2$-boundedness
of $\phi^n$, there
 exists a constant $C$ such that
$$\|z^n -z\|
\le C\left(\sup_{s\in[0,T]}\| h^n_s\|_H+
 \sup_{s\in[0,T]}\|
h^n_s\|_H^2+ \left( \int_0^T \|h^n_s\|_V^\alpha\d
 s\right)^{\frac{1}{\alpha}}\right). $$
Since $h^n\rightarrow 0$ in $C([0,T]; H)\cap L^\alpha([0,T]; V)$, we
know
$z^n\rightarrow z \  \text{strongly \ in}\
C([0,T]; H)\cap L^\alpha([0,T]; V)$
 as $n\rightarrow\infty$.

 Since Lemma
 \ref{L2.1} shows that the convergence of the corresponding solution
  $z^{\phi}$ is uniformly on $S_N$ w.r.t. the approximation on $\phi$,  the
 conclusion on the
 case  $\phi^n, \phi\in L^2([0,T]; U)$  can de derived by the proof above   and
the standard $3\vv$-argument.

\textbf{Step 2}: Now we prove the conclusion for general $B$ without assuming
 $(A5)$. Denote $z_{t,n}^\phi$ the solution of the
following equation
$$\ff{\d z_{t,n}^\phi}{\d t} = A(t,z_{t,n}^\phi )  +P_n
 B(t,z_{t,n}^\phi)\phi_t,\
\ \ z_{0,n}^\phi= x, $$ where $P_n$ is the standard projection(see $(A4)$ and Section 4
 for details). By using the same argument in Lemma \ref{L2.1} we can
 prove
\beq
\begin{split}\label{finite dim estimate}
  & \|z^\phi_{n} -z^\phi\|^2+\delta\int_0^T\|z^\phi_{s,n}
 -z_s^\phi\|_V^\alpha\d s\\
& \le \exp\left\lbrace \int_0^T(K+2\|\phi_s\|_U^2)\d s\right\rbrace
\int_0^T\|(I-P_n)B(s,z_s^\phi)\|_2^2\d s
\end{split}\end{equation}
Since $B(\cdot,\cdot)$ are Hilbert-Schmidt (hence compact)
operators, then by the dominated convergence theorem we know
$$ \int_0^T\|(I-P_n)B(s,z_s^\phi)\|_2^2\d s\rightarrow 0\ \text{as}\ n\rightarrow\infty.$$
Hence $z^\phi_{n}\rightarrow z^\phi$ in $C([0,T]; H)\cap
L^\alpha([0,T]; V)$ as $n\rightarrow\infty$.
 Moreover, this convergence is uniformly (w.r.t $\phi$)
on bounded set of $L^2([0,T]; U)$, which follows from (\ref{finite
dim
 estimate}) and (\ref{bound of B}).

 Note that $P_nB$ satisfy $(A5)$, by combining with \textbf{Step
1}
 and standard $3\varepsilon$-argument we can conclude that
$z^n\rightarrow z$ strongly in $C([0,T]; H)\cap L^\alpha([0,T]; V)$
 for general $B$.
Hence the proof is complete.
\end{proof}


\beg{lem} \label{L2.3} Assume $(A1)-(A3)$ and $(A5)$ hold. Let
$\{v^\varepsilon\}_{ \varepsilon>0}\subset \mathcal{A}_N$ for some
 $N<\infty$. Assume $v^\varepsilon$ converge to $v$ in distribution as
 $S_N$-valued random elements, then
$$g^\varepsilon\left(W_\cdot+\frac{1}{\varepsilon}
\int_0^\cdot v^\varepsilon_s\d s\right)\rightarrow
g^0\left(\int_0^\cdot
 v_s\d s \right) $$
in distribution as $\vv\rightarrow 0$.
\end{lem}
\begin{proof} By the Girsanov theorem and uniqueness of solution to
 (\ref{1.3}), it's easy to see that
 $X^{\varepsilon}:=g^\varepsilon\left(W_\cdot+\frac{1}{\varepsilon}
\int_0^\cdot v^\varepsilon_s\d s\right)$ (the abuse of notation here
is for
 simplicity) is the unique solution of the following equation
\beq\label{app}
  \d X_t^\varepsilon =
\left(A(t,X_t^\varepsilon)+B(t,X_t^\varepsilon)v_t^\varepsilon\right)\d
t+ \varepsilon B(t,X_t^\varepsilon)\d W_t, \
 X_0^\vv=x.
\end{equation}
Now we only need to show $X^{\varepsilon}\rightarrow z^v$ in
 distribution as $\varepsilon\rightarrow0$.
We may assume $\vv\leq\frac{1}{2}$, by using the It\^o formula,
Young's
 inequality and $(A2)$ we have
\beq\label{compare}\begin{split}
 \d  \|X_t^\varepsilon -z_t^v\|_H^2
= & 2{ }_{V^*}\< A(t,X_t^\vv )
-A(t,z_t^v), X_t^\vv-z_t^v\>_V\d t\\
& +2\< X_t^\vv-z_t^v,
 (B(t,X_t^\vv)-B(t,z_t^v))v^\vv_t+B(t,z_t^v)(v^\vv_t-v_t)\>_H \d t\\
&  +\varepsilon^2\|B(t,X_t^\varepsilon)\|_{2}^2 \d t
+2\vv\<X_t^\varepsilon -z_t^v,   B(t,X_t^\varepsilon)\d W_t\>_H\\
\le & \left( 2{ }_{V^*}\< A(t,X_t^\vv ) -A(t,z_t^v),
 X_t^\vv-z_t^v\>_V+\|B(t,X_t^\vv)-B(t,z_t^v)\|_{2}^2\right)\d t\\
& +2\|v^\vv_t\|_U^2\|X_t^\vv-z_t^v\|_H^2\d t +2\< X_t^\vv-z_t^v,
 B(t,z_t^v)(v^\vv_t-v_t)\>_H  \d t\\
&  +2\varepsilon^2\|B(t,z_t^v)\|_{2}^2 \d t
+2\vv\<X_t^\varepsilon -z_t^v,   B(t,X_t^\varepsilon)\d W_t\>_H\\
\le & \left[-\delta\|X_t^\varepsilon -z_t^v\|_V^\alpha +
C(1+\|v_t^\vv\|_U^2)\|X_t^\varepsilon
-z_t^v\|_H^2+ 2\varepsilon^2\|B(t,z_t^v)\|_{2}^2 \right] \d t\\
& +2\< X_t^\vv-z_t^v, B(t,z_t^v)(v^\vv_t-v_t)\>_H  \d t
 +2\vv\<X_t^\varepsilon -z_t^v,  B(t,X_t^\varepsilon)\d W_t\>_H.
\end{split}\end{equation}
Similarly we define
$$h_t^\vv=\int_0^t B(s,z_s^v)(v^\vv_s-v_s)\d s,$$
then we know that $h^\vv\rightarrow 0$ in distribution as $C([0,T];
 V)$-valued random element, consequently
 also in $C([0,T];
 H)\cap L^\alpha([0,T]; V)$.
Note that
$$2\<X_t^\vv-z_t^v, h^\vv_t\>_H=\|X_t^\vv-z_t^v+
 h^\vv_t\|_H^2-\|X_t^\vv-z_t^v\|_H^2-\|h^\vv_t\|_H^2.$$
By using the It\^o formula for corresponding square norm we can
 derive that
\begin{equation}\begin{split}\label{e3}
& \int_0^t \<X_s^\vv-z_s^v, B(s,z_s^v)(v^\vv_s-v_s)\>_H\d s\\
 =&
\<X_t^\vv-z_t^v, h^\vv_t\>_H-\int_0^t{
 }_{V^*}\<A(s,X_s^\vv)-A(s,z_s^v), h^\vv_s\>_V\d s\\
& -\int_0^t \<B(s,X_s^\vv)v_s^\vv-B(s,z_s^v)v_s, h^\vv_s\>_H\d s
 -\vv\int_0^t \<B(s,X_s^\vv)\d W_s, h^\vv_s\>_H.
\end{split}
\end{equation}
By using the same argument as in (\ref{e2}) we obtain
\begin{equation}\label{e4}
\beg{split}
& \int_0^t \<X_s^\vv-z_s^v, B(s,z_s^v)(v^\vv_s-v_s)\>_H\d s\\
\le & \frac{1}{4}\|X_t^\vv-z_t^v\|_H^2+ \sup_{s\in[0,t]}\|
h^\vv_s\|_H^2-\vv\int_0^t \<B(s,X_s^\vv)\d W_s,
 h^\vv_s\>_H\\
& + C\left(\int_0^t \left(
 1+\|z_s^v\|_V^{\alpha}+\|X_s^\vv\|_V^{\alpha}\right)
\d s\right)^{\frac{\alpha-1}{\alpha}}\cdot
\left( \int_0^t \|h^\vv_s\|_V^\alpha\d s\right) ^{\frac{1}{\alpha}}\\
& +C\sup_{s\in[0,t]}\| h^\vv_s\|_H \left\lbrace \left( \int_0^t
\|B(s,X_s^\vv)\|_2^2\d s\right)^{1/2} + \left( \int_0^t
\|B(s,z_s^v)\|_2^2\d s\right)^{1/2} \right\rbrace.
\end{split}
\end{equation}
Hence from (\ref{compare})-(\ref{e4}) we have
\begin{equation}\beg{split}
& \|X_t^\varepsilon -z_t^v\|_H^2+\delta\int_0^t\|X_t^\varepsilon
 -z_t^v\|_V^\alpha\d s\\
& \le c_1\int_0^t(1+\|v_s^\vv\|_U^2)\|X_t^\varepsilon-z_t^v\|_H^2\d
 s+c_2(\vv^2+\sup_{s\in[0,t]}\|
h^\vv_s\|_H^2)\\
& + c_3\left(1+\int_0^t \|X_s^\vv\|_V^{\alpha} \d
s\right)^{\frac{\alpha-1}{\alpha}}\cdot
\left( \int_0^t \|h^\vv_s\|_V^\alpha\d s\right) ^{\frac{1}{\alpha}}\\
& +c_4\sup_{s\in[0,t]}\| h^\vv_s\|_H \left\lbrace 1+\left(\int_0^t
\|X_s^\vv\|_H^2 \d s\right)^{1/2}
\right\rbrace  \\
& +
 4\vv\int_0^t\<X_s^\varepsilon -z_s^v-h^\vv_s,  B(s,X_s^\varepsilon)\d
 W_s\>_H,
\end{split}\end{equation}
where we used the estimate (see (\ref{bound})-(\ref{bound of B}))
that
 there exists constant $C$ such that
$$\int_0^T\|B(s,z_s^v)\|_{2}^2\d s+\int_0^T
\|z_s^v\|_V^{\alpha}\d s\le C,\ \ a.s. $$
By applying the Gronwall lemma
we have \beq\label{e5}\beg{split} &
\sup_{s\in[0,t]}\|X^\varepsilon_s -z^v_s\|_H^2
 +\delta\int_0^t\|X_s^\varepsilon-z_s^v \|_V^{\alpha}\d s\\
& \le  C\bigg[\varepsilon^2+\sup_{s\in[0,t]}\| h^\vv_s\|_H^2
  +\left(1+\int_0^t \|X_s^\vv\|_V^{\alpha} \d
s\right)^{\frac{\alpha-1}{\alpha}}
\left( \int_0^t \|h^\vv_s\|_V^\alpha\d s\right) ^{\frac{1}{\alpha}}\\
& +\sup_{s\in[0,t]}\| h^\vv_s\|_H \left\lbrace 1+\left(\int_0^t
\|X_s^\vv\|_H^2\d
 s\right)^{1/2}
\right\rbrace
 +\sup_{u\in[0,t]}\left|\vv\int_0^u\<X_s^\varepsilon
-z_s^v-h^\vv_s,
  B(s,X_s^\varepsilon)\d W_s\>_H \right| \bigg]
\end{split}\end{equation}
Define the stopping time
$$ \tau^{M,\vv}=\inf \left\lbrace t\le T:
 \sup_{s\in[0,t]}\|X^\vv_s\|_H^2+\int_0^t\|X^\vv_s\|_V^\alpha\d s>M
 \right\rbrace. $$
By the Burkh\"older-Davis-Gundy inequality one has
\beq\label{e6}\beg{split} &
\vv\textbf{E}\sup_{t\in[0,\tau^{M,\vv}]}\left|\int_0^t\<X_s^\varepsilon
 -z_s^v-h^\vv_s,  B(s,X_s^\varepsilon)\d W_s\>_H \right|\\
\le & 3\vv\textbf{E}\left\lbrace
\int_0^{\tau^{M,\vv}}\|X_s^\varepsilon
 -z_s^v-h^\vv_s\|_H^2 \|B(s,X_s^\varepsilon)\|_2^2\d s\right\rbrace
 ^{1/2}\\
\le & 3\vv\textbf{E}\left\lbrace
 \sup_{s\in[0,\tau^{M,\vv}]}\|X_s^\varepsilon
 -z_s^v-h^\vv_s\|_H^2+C\int_0^{\tau^{M,\vv}}
\left( 1+\|X_s^\varepsilon\|_H^2+\|X_s^\varepsilon\|_V^\alpha\right) \d
 s\right\rbrace \\
\le & C\vv\rightarrow 0 \  (\vv\rightarrow 0).
\end{split}\end{equation}

By using the similar argument in
(\ref{compare})
 we have
$$\d \|X_t^\vv\|_H^2\le
-\frac{\delta}{2}\|X_t^\vv\|_V^\alpha\d t+C(1+\|X_t^\vv\|_H^2+\|v_t^\vv\|_U^2\|X_t^\vv\|_H^2)\d
t +2\vv\<X_t^\vv, B(t,X_t^\vv)\d W_t\>_H,$$ where $C$ is a constant.
Repeat the same argument in \cite[Theorem 3.10]{KR}  we can prove
$$\sup_{\vv\in[0,1)}\textbf{E}\left\lbrace
\sup_{t\in[0,T]}\|X_t^\varepsilon\|_H^2+\int_0^T\|X_t^\varepsilon\|_V^\alpha\d
t\right\rbrace <\infty . $$
Hence there exists a  suitable constant $C$ such that
\beq\label{e8}
\liminf_{\vv\rightarrow0}\textbf{P}\{\tau^{M,\vv}=T\}\geq
1-\frac{C}{M}.
\end{equation}
Recall that $h^\vv\rightarrow 0$ in distribution in $C([0,T]; H)\cap
 L^\alpha([0,T]; V)$, combining with (\ref{e5})-(\ref{e8}) one can
 conclude
$$\sup_{t\in[0,T]}\|X^\vv_t -z^v_t\|_H^2+\int_0^T \|X_t^\vv-z_t^v
 \|_V^{\alpha}\d t
 \rightarrow0\ (\vv\rightarrow 0) $$
in distribution. Hence the proof is complete.
\end{proof}
\beg{rem}\label{rem2} According to Lemma \ref{L1}, Lemma \ref{L2.2}
and Lemma \ref{L2.3}, we know that $\{X^\vv\}$ satisfy LDP provided $(A1)-(A3)$
and $(A5)$ hold. By using some approximation argument, we can
  replace  $(A5)$ by the weaker assumption $(A4)$.
\end{rem}

\section{Replace $(A5)$ by $(A4)$}

 For any fixed $n\ge 1$, let
$H_n\subseteq V$ compact and $P_n: H\to H_n$ be the orthogonal
projection.  Let $X_{t}^{\vv,n}$ be the solution of
\beq\label{finite} \d X_{t}^{\vv,n}= A(t,X_t^{\vv,n})\d t +\vv P_n
 B(t,X_t^{\vv,n})\d W_t,\ \ \ X_0^{\vv, n}=x.\end{equation}
Since $P_nB$ satisfy $(A5)$, according to the Section 3(Remark
\ref{rem2}) we know $\{X^{\vv,n}\}$ satisfy the LDP provided
$(A1)-(A3)$. Now we prove that $\{X^{\vv,n}\}$ are the exponential
good approximation to
 $\{X^{\vv}\}$ if the following assumption hold.
\begin{enumerate}
   \item [$(A4^\prime)$]
$$a_n:=\sup_{(t,v)\in
 [0,T]\times V}\|P_nB(t,v)-B(t,v)\|_2^2\rightarrow 0\ (n\rightarrow\infty).$$
\end{enumerate}

\beg{lem}\label{finite-dim approximation}If $(A1)-(A3)$ and $(A4^\prime)$
 hold, then $\forall
\sigma>0$
 \beq\label{W2} \limsup_{n\to\infty}\limsup_{\vv\to
0}\vv^2\log \mathbf{P} \left(\rho(X^{\vv},
X^{\vv,n})>\sigma\right)=-\infty,
\end{equation}
where $\rho$ is the metric on $C([0,T];H)\cap L^\alpha([0,T];V)$
defined in $(\ref{metric})$.
\end{lem}

\begin{proof} For $\vv<\frac{1}{2}$,
 by using the It\^o formula and $(A2)$
 we have
$$
\aligned
& \d \|X_t^\vv -X_t^{\vv,n}\|_H^2 \\
 =& \left( 2{ }_{V^*}\< A(t,X_t^\vv)-A(t,X_t^{\vv,n}),
X_t^\vv-X_t^{\vv,n}\>_V+\varepsilon^2\|B(t,X_t^\varepsilon)-P_nB(t,X_t^{\vv,n})\|_{2}^2
 \right) \d t\\
&   +2\vv\<X_t^\varepsilon -X_t^{\vv,n},
 (B(t,X_t^\varepsilon)-P_nB(t,X_t^{\vv,n}))\d W_t\>_H,
\endaligned $$
where  $C$ is a constant.  Define
$$
\|X_t^\vv-X_t^{\vv,n}\|=\|X_t^\vv-X_t^{\vv,n}\|_H^2+\delta\int_0^t\|X_s^\vv
-X_s^{\vv,n}\|_V^\alpha\d s.$$
Note that
$$ M_t^{(n)}:=\int_0^t \<X_s^\vv-X_s^{\vv,n},\left(B(s,X_{s}^{\vv})-P_n
B(s,X_{s}^{\vv,n})\right)\d W_t\>_H$$
 is a local martingale and its quadratic variation process satisfies
$$
\d \<M^{(n)}\>_t\le
 2\|X_t^\vv-X_t^{\vv,n}\|_H^2(\|B(t,X_t^\vv)-B(t,X_t^{\vv,n})\|_2^2 +a_n) \d t.
$$
Let $\varphi_\theta(y)=(a_n+y)^\theta$ for some $\theta\le \frac{1}{4\varepsilon^2}$, then by $(A2)$
\beq\beg{split}
& \d \varphi_\theta(\|X_t^\vv-X_t^{\vv,n}\|)\\
\le&
 \theta(a_n+\|X_t^\vv-X_t^{\vv,n}\|)^{\theta-1}
\left(\d \|X_t^\vv -X_t^{\vv,n}\|_H^2+\delta \|X_t^\vv
-X_t^{\vv,n}\|_V^\alpha\d t
 \right) \\
& +4\vv^2\theta(\theta-1)(a_n+\|X_t^\vv-X_t^{\vv,n}\|)^{\theta-2}
\|X_t^\vv-X_t^{\vv,n}\|_H^2\left(\|B(t,X_t^\vv)-B(t,X_t^{\vv,n})\|_2^2
 +a_n\right) \d t\\
\le &
 C\theta \varphi_\theta\left(\|X_t^\vv-X_t^{\vv,n}\|\right)\d t+\d
 \beta_t
\end{split}\end{equation}
where $C$ is a constant and $\beta_t$ is a local martingale. By standard localization
 argument we may assume $\beta_t$ is a martingale  for simplicity.
Let $\theta=\frac{1}{4\vv^2}$ we know
$$N_t:=\exp\left[-\frac{C}{4\vv^2}t \right]\varphi_{\frac{1}{4\vv^2}}\left(\|X_t^\vv-X_t^{\vv,n}\|\right)$$
is a  supermartingale. Hence we have
$$\aligned
& \mathbf{P}\left( \rho(X^{\vv},
X^{\vv,n})>2\sigma \right) \\
 \le & \mathbf{P}\left(
 \sup_{t\in[0,T]}\|X_t^{\vv}-X_t^{\vv,n}\|_H>\sigma\right)
+\mathbf{P}\left(\int_0^T\|X_t^{\vv}-X_t^{\vv,n}\|_V^\alpha\d t>\sigma^\alpha\right)\\
\le & \mathbf{P}
\left(\sup_{t\in[0,T]}N_t>\exp\left[-\frac{C}{4\vv^2}T\right](\sigma^2+a_n)^{\frac{1}{4\vv^2}}\right)
+\mathbf{P}
\left(\sup_{t\in[0,T]}N_t>\exp\left[-\frac{C}{4\vv^2}T\right](\delta\sigma^\alpha+a_n)^{\frac{1}{4\vv^2}}\right)\\
\le&
\exp\left[\frac{C}{4\vv^2}T\right](\sigma^2+a_n)^{-\frac{1}{4\vv^2}}\textbf{E}
 N_0
+\exp\left[\frac{C}{4\vv^2}T\right](\delta\sigma^\alpha+a_n)^{-\frac{1}{4\vv^2}}\textbf{E}
  N_0\\
=& \exp\left[\frac{C}{4\vv^2}T\right]\left\{ \left(
 \frac{a_n}{\sigma^2+a_n}\right)^{\frac{1}{4\vv^2}}+
\left(\frac{a_n}{\delta\sigma^\alpha+a_n}\right)^{\frac{1}{4\vv^2}} \right\}.
\endaligned$$
This implies that
$$\aligned
& \limsup_{\vv\to 0}\vv^2\log
\mathbf{P}\left(\rho(X^{\vv}, X^{\vv,n})>2\sigma\right)\\
& \le \frac{CT}{4}+ \max\left\lbrace \log \frac{a_n}{\sigma^2+a_n}, \log
 \frac{a_n}{\delta\sigma^\alpha+a_n}\right\rbrace.
\endaligned $$
Since $(A4^\prime)$ says $a_n\rightarrow0$ as $n\rightarrow\infty$, $(\ref{W2})$ hold and the proof is
complete.

\end{proof}

\begin{cor}\label{cor1}
If $(A1)-(A3)$ and $(A4^\prime)$ hold, then $\{X^\vv\}$
  satisfy the LDP in $C([0,T]; H)\cap
L^\alpha([0,T]; V)$ with rate function $(\ref{rate})$.
\end{cor}
\begin{proof}
According to \cite[Theorem 2.1]{Wu04} and section 3 one can conclude
$\{X^\vv\}$
  satisfy the LDP with the following rate function
$$\tilde{I}(f):=\sup_{r>0}\liminf_{n\rightarrow\infty}\inf_{g\in
 S_r(f)}I^n(g)=
\sup_{r>0}\limsup_{n\rightarrow\infty}\inf_{g\in S_r(f)}I^n(g).$$
where $S_r(f)$ is the closed ball  in $C([0,T]; H)\cap
L^\alpha([0,T]; V)$ centered at $f$ with
 radius $r$ and  $I^n$ is
 given by
\beq\label{rate 3} I^n(z):=  \inf \left\{\ff 1 2
\int_0^T\|\phi_s\|_{U}^2\d
 s:\
z=z^{n,\phi} ,\ \phi\in L^2([0,T], U)\right\},
\end{equation}
where $z^{n,\phi}$ is the unique solution of following equation
$$\frac{\d z_{t}^n}{\d t} = A(t,z_{t}^n)+ P_nB(t,z_{t}^n)\phi_t, \
 z_{0}^n=x. $$
Now we only need to prove $\tilde{I}=I$, i.e.
$$I(f)=\sup_{r>0}\liminf_{n\rightarrow\infty}\inf_{g\in
 S_r(f)}I^n(g).$$
We will first show that for any $r>0$
$$I(f)\ge \liminf_{n\rightarrow\infty}\inf_{g\in S_r(f)}I^n(g).$$
We assume $I(f)<\infty$, then by Lemma \ref{L2.2} there exists
$\phi$
 such that
$$f=z^\phi\ \ \text{and}\ \ I(f)=\frac{1}{2}\int_0^T\|\phi_s\|_{U}^2\d
 s.$$
Since $z^{n,\phi}\rightarrow z^\phi$, for $n$ large enough we have
 $$f_n:=z^{n,\phi}\in S_r(f).$$
Notice $I^n(f_n)\le\frac{1}{2}\int_0^T\|\phi_s\|_{U}^2\d s$, hence
we
 have
$$\liminf_{n\rightarrow\infty}\inf_{g\in S_r(f)}I^n(g)\le
\liminf_{n\rightarrow\infty}I^n(f_n)\le I(f) .$$ Since $r$ is
arbitrary we have proved the lower bound
$$I(f)\ge\sup_{r>0}\liminf_{n\rightarrow\infty}\inf_{g\in
 S_r(f)}I^n(g).$$
For the upper bound we can proceed as in finite dimensional case in
 \cite[Lemma 4.6]{St} to show
$$\limsup_{n\rightarrow\infty}\inf_{g\in S_r(f)}I^n(g)\ge \inf_{g\in
 S_r(f)}I(g) $$
Hence we have
$$\sup_{r>0}\limsup_{n\rightarrow\infty}\inf_{g\in S_r(f)}I^n(g)\ge
 \sup_{r>0}\inf_{g\in S_r(f)}I(g)\ge I(f).$$
Hence the proof is complete.
\end{proof}

In order to replace the assumption $(A4^\prime)$ by $(A4)$, we need to
use some truncation techniques
(cf. \cite{St,C}).

\beg{lem}\label{L4}Assume $(A1)-(A4)$ hold, then \beq\label{W3}
\lim_{R\to\infty}\limsup_{\vv\to 0}\vv^2\log \mathbf P
(\sup_{t\in[0,T]}\|X_t^{\vv}\|_H^2+\frac{\delta}{2}\int_0^T\|X_t^\vv\|_V^\alpha\d
 t>R)=-\infty .
\end{equation}
\end{lem}

\begin{proof}
 By using the It\^o formula  we
 have
$$
\aligned \d \|X_t^\vv\|_H^2 & = \left( 2{ }_{V^*}\< A(t,X_t^\vv),
 X_t^\vv\>_V+\varepsilon^2\|B(t,X_t^\varepsilon)\|_{2}^2 \right) \d t
+2\vv\<X_t^\varepsilon,   (B(t,X_t^\varepsilon)\d W_t\>_H.
\endaligned $$
Note that  $M_t^{(n)}:=\int_0^t \<X_s^\vv, B(s,X_{s}^{\vv})\d W_s\>_H$
is a local martingale and
$$
d\<M^{(n)}\>_t\le \|X_t^\vv\|_H^2\|B(t,X_t^\vv)\|_2^2 \d t.
$$
Define
$$\|X_t^\vv\|:=
 \|X_t^\vv\|_H^2+\frac{\delta}{2}\int_0^t\|X_s^\vv\|_V^\alpha\d s, \ \
\varphi_\theta(y)=(1+y)^\theta,\ \theta>0,$$
then  for $\theta\le\frac{1}{2\varepsilon^2}$ by $(A2)$ and $(A3)$ we have
\beq\label{e}\beg{split} \d \varphi_\theta(\|X_t^\vv\|)\le&
\theta(1+\|X_t^\vv\|)^{\theta-1}
\left(\d\|X_t^\vv\|_H^2+\frac{\delta}{2}\|X_t^\varepsilon\|_V^\alpha \d t\right) \\
& +2\vv^2\theta(\theta-1)(1+\|X_t^\vv\|)^{\theta-2}
\|X_t^\vv\|_H^2 \|B(t,X_t^\vv)\|_2^2 \d t\\
\le & C\theta\varphi_\theta(\|X_t^\vv\|)\d t+\d
\beta_t
\end{split}\end{equation}
where $\beta_t$ is a local martingale. We also omit the
standard
 localization procedure here. Let $\theta=\frac{1}{2\vv^2}$ we know
$$N_t:=\exp\left[-\frac{C}{2\vv^2}t\right]\varphi_{\frac{1}{2\vv^2}}\left(\|X_t^\vv\|\right)$$
is a  supermartingale. Hence we have
$$\aligned
&
 \mathbf{P}\left(\sup_{t\in[0,T]}\|X_t^{\vv}\|_H^2+\frac{\delta}{2}\int_0^T\|X_t^\vv\|_V^\alpha\d t>R\right)\\
\le & \mathbf{P}
\left(\sup_{t\in[0,T]}N_t>\exp\left[-\frac{C}{2\vv^2}T\right](1+R)^{\frac{1}{2\vv^2}}\right)\\
\le & \exp\left[\frac{C}{2\vv^2}T\right](1+R)^{-\frac{1}{2\vv^2}}\textbf{E} N_0\\
= & \exp\left[\frac{C}{2\vv^2}T\right]\left( \frac{1}{1+R}\right)
 ^{\frac{1}{2\vv^2}}.
\endaligned$$
This implies that
$$\limsup_{\vv\to 0}\vv^2\log
\mathbf{P}\left(\sup_{t\in[0,T]}\|X_t^{\vv}\|>R\right)\le \frac{1}{2}\log
 \frac{1}{1+R}+\frac{CT}{2}.  $$
Therefore, (\ref{W3}) hold.
\end{proof}

After all these preparations, now we can finish the proof of Theorem
\ref{T1.1}.

\noindent\textbf{Proof of Theorem \ref{T1.1}}: The proof here is a
slight
 modification of \cite[Theorem 4.13]{St}.
Define $\xi:V\rightarrow[0,1]$ be a $C_0^\infty$-function such that
$$\xi(v):= \beg{cases} 0, &\text{if}\  \|v\|_H>2,\\
 1, &\text{if}\  \|v\|_H\le 1.\end{cases}$$
Let $\xi_N(v)=\xi(\frac{v}{N})$ and
$$B_N(t,v)=\xi_N(v)B(t,v)+(1-\xi_N(v))B(t,0).$$
Consider the mollified problem for equation (\ref{1.3}):
\beq\label{mollified equation}
  \d X_{t,N}^\vv = A(t,X_{t,N}^\vv)\d t+\vv B_N(t,X_{t,N}^\vv)\d W_t, \
 X_0=x.
\end{equation}
It's easily to see that $A, B_N$ satisfy $(A1)-(A3)$ and $(A4^\prime)$, since
in this case $(A4)$ implies that for $B_N$
$$a_n=\max\left\{\sup_{(t,v)\in[0,T]\times
S_{2N}}\|(I-P_n)B(t,v)\|_{2}^2,\
\sup_{t\in[0,T]}\|(I-P_n)B(t,0)\|_{2}^2\right\} \rightarrow 0
(n\rightarrow \infty).
$$
Hence by Corollary \ref{cor1}   we know $\{X_{N}^\vv\}_{\vv>0}$
 satisfy large deviation principle on $C([0,T];H)\cap
 L^\alpha([0,T]; V)$ with the following mollified rate function
\beq\label{rate 2} I_N(z):=\inf
\left\{\frac{1}{2}\int_0^T\|\phi_s\|_{U}^2\d
  s:
\ z=z_N^{\phi} ,\ \phi\in L^2([0,T], U)\right\},
\end{equation}
where $z_N^{\phi}$ is the unique solution of following equation
$$\frac{\d z_{t,N}}{\d t} = A(t,z_{t,N})+ B_N(t,z_{t,N})\phi_t, \
 z_{0,N}=x.  $$
Let $N\rightarrow\infty$, then the LDP for $\{X^\vv\}$ can be
 derived as in the finite dimensional case.

 According to Lemma \ref{L2.2}, $I$ defined in (\ref{rate})
is
 a (good) rate function.
Note  $I_N(z)=I(z)$ for any $z\in C([0,T];H)\cap L^\alpha([0,T];
 V)$ satisfy
 $$\|z\|_T:=\sup_{t\in[0,T]}\|z_t\|_H\le N.$$
  We now first
 show
 that for any open set $G\subseteq C([0,T];H)\cap L^\alpha([0,T];
 V)$
$$\liminf_{\vv\rightarrow0}\vv^2\log \mathbf{P}\left(X^\vv\in
 G\right)\ge -\inf_{z\in G}I(z).$$
Obviously, we only need to prove that for all  $\overline{z}\in G$
with
 $\overline{z}_0=x$
$$\liminf_{\vv\rightarrow0}\vv^2\log \mathbf{P}\left(X^\vv\in
 G\right)\ge -I(\overline{z}).$$
Choose $R>0$ such that $\|\overline{z}\|_T<R$ and set
$$N_R=\{z\in C([0,T];H)\cap L^\alpha([0,T]; V):\|z\|_T<R\}.$$
Then we have
$$\aligned
 \liminf_{\vv\rightarrow0}\vv^2\log \mathbf{P}\left(X^\vv\in
 G\right)&\ge
\liminf_{\vv\rightarrow0}\vv^2\log \mathbf{P}\left(X^\vv\in G\cap
 N_R\right)\\
& =  \liminf_{\vv\rightarrow0}\vv^2\log \mathbf{P}\left(X^\vv_N\in
 G\cap N_R\right)\\
& \ge -\inf_{z\in G\cap N_R}I_N(z)\\
& \ge -I(\overline{z}).
\endaligned$$
Finally, given a closed set $F$ and an $L<\infty$, by Lemma \ref{L4}
 there exists $R$ such that
$$\aligned
 \limsup_{\vv\rightarrow0}\vv^2\log \mathbf{P}\left(X^\vv\in
 F\right)&\le
\limsup_{\vv\rightarrow0}\vv^2\log\left(\mathbf{P}(X^\vv\in F\cap
 \overline{N_R})+\mathbf{P}(X^\vv\in N_R^c)\right)\\
& \le (-\inf_{z\in F\cap \overline{N_R}}I_N(z))\vee (-L)\\
& \le -\left[ \inf_{z\in F}I(z)\wedge L\right].
\endaligned$$
Let $L\rightarrow\infty$, we obtain
$$\limsup_{\vv\rightarrow0}\vv^2\log \mathbf{P}\left(X^\vv\in
 F\right)\le -\inf_{z\in F}I(z). $$
Now the proof of Theorem \ref{T1.1} is complete.  \qed

\section{Examples}
Now we can apply the main results to many stochastic evolution equations
as applications. As a preparation  we prove the following lemma
first.

\beg{lem}\label{L6.1} Let $(E, \langle\cdot,\cdot\rangle, \|\cdot\|)$ is a
Hilbert
 space,
then for any $r\geq0$ we have \beq \label{6.1}
\langle\|a\|^ra-\|b\|^rb, a-b\rangle \geq2^{-r}\|a-b\|^{r+2},\  a,b\in
E. \end{equation}
 \beq\label{6.2} \|\|a\|^{r-1}a-\|b\|^{r-1}b\| \leq \max\{r,
1\}\|a-b\|(\|a\|^{r-1}+\|b\|^{r-1}),\  a,b \in E.
\end{equation}
If  $0<r<1$, then there exists a constant $C>0$ such that
 \beq\label{6.3} ||a|^{r-1}a-|b|^{r-1}b| \leq C |a-b|^{r},\  a,b \in
 \mathbb{R}.
\end{equation}

 \end{lem}
\begin{proof}
(i) By the symmetry of (\ref{6.1})  we may assume $\|a\|\geq\|b\|$. Then
 \ce
&&\langle\|a\|^ra-\|b\|^rb,
a-b\rangle\\
&=&\|b\|^r\|a-b\|^2+(\|a\|^r-\|b\|^r)\langle a, a-b\rangle\\
&=&\|b\|^r\|a-b\|^2+(\|a\|^r-\|b\|^r)\cdot\frac{1}{2}(\|a\|^2+\|a-b\|^2-\|b\|^2)\\
&\geq&\|b\|^r\|a-b\|^2+\frac{1}{2}(\|a\|^r-\|b\|^r)\|a-b\|^2\\
&=&\frac{1}{2}(\|a\|^r+\|b\|^r)\|a-b\|^2\\
&\geq&2^{-r}\|a-b\|^{r+2}, \de
since $\|a-b\|^r\leq2^{r-1}(\|a\|^r+\|b\|^r)$.

 (ii) The proof of (\ref{6.2}) and (\ref{6.3}) is similar.
\end{proof}

The first example is to obtain the  LDP for a class of
reaction-diffusion type SPDEs within the variational framework, which improve the main result in \cite{C}.
\beg{exa} (Stochastic reaction-diffusion equations)\\
 Let $\Lambda$ be an open bounded domain in $\mathbb{R}^d$
with smooth boundary and $L$ be a negative definite self-adjoint
operator on $H:=L^2(\Lambda)$. Suppose
$$V:=\D(\sqrt{-L}),\ \    \|v\|_V:=\|\sqrt{-L}v\|_H.$$
is a Banach space such that $V\subseteq H$ is dense and compact, and $L$ can
be extended as a continuous operator from  $V$ to it's dual space
$V^*$. Consider the following  semilinear stochastic equation
 \beq\label{rd}
 \d X_t^\vv=(LX_t^\vv+F(t,X_t^\vv))\d t+\vv B(t,X_t^\vv)\d W_t,
\ X_0^\vv=x\in H,
  \end{equation}
where $W_t$ is a cylindrical Wiener process on another separable
Hilbert
 space $U$ and
$$F:[0,T]\times
V\rightarrow V^*,\ \  B:[0,T]\times V\rightarrow L_2(U;V).$$ If $F$
and $B$
 satisfy the following conditions:
 \begin{equation}\begin{split}\label{assumption for rd}
2{ }_{V^*}\langle F(t,u)-F(t,v), u-v\rangle_{V}&+\|B(t,u)-B(t,v)\|_{2}^2  \le
 C\|u-v\|_H^2,\\
\|F(t,v)\|_{V^*} \le
 C(1+\|v\|_V),&\ \ \|B(t,v)\|_{2} \le
 C(1+\|v\|_H^\gamma), \ u,v\in V.
\end{split}\end{equation}
where $C,\gamma>0$ are constants,
 then $\{X^\vv\}$ satisfy the large deviation principle on
 $C([0,T];H)\cap L^2([0,T];V)$.

 \end{exa}
 \beg{proof} From the assumptions (\ref{assumption for rd}), it's easy to show that $(A1)-(A4)$
 hold for $\alpha=2$.
Hence the conclusion follows from Theorem \ref{T1.1}.
 \end{proof}

  \beg{rem}
 (i) We can simply take $L$ as the Laplace operator with Dirichlet
 boundary
 condition and  $F(t,X_t)=-|X_t|^{p-2}X_t (1\leq p\leq2)$ as
 a concrete example.

(ii) Compare with the result in \cite[Theorem 4.2]{C}(only time
homogeneous case), the author in
 \cite{C}
need to assume $F$ is local Lipschitz and have more restricted range
conditions:
$$F:[0,T]\times V\rightarrow H.$$
In our example we can allow $F$ is monotone and take values in
$V^*$. Another difference is we also drop the non-degenerated
condition $(A.4)$ on $B$ in \cite{C}.

(iii) Note here one can also take $B: V\rightarrow L_2(U;H)$ with
locally compact range, which seems not allowed in \cite[Theorem
4.2]{C}.
 \end{rem}

The second example is  stochastic  porous media equations,
which have been studied  intensively in recent years, see
e.g.\cite{DRRW,RRW,RWW,Wang}. We use the same framework
 as in \cite{RWW,Wang}.
\beg{exa}\label{E6.2} (Stochastic porous media equations)\\
 Let $(E,\mathcal{M},{\bf m})$ be a separable probability space and $(L,\D(L))$ a
negative definite self-adjoint linear operator on $(L^2({\bf
m}),\<\cdot ,\cdot \>)$ with spectrum contained in
$(-\infty,-\ll_0]$ for some $\ll_0>0$.  Then the embedding
$$H^1:=\D(\sqrt{-L})\subseteq L^2(\bf m)$$
 is dense and continuous.
Define $H$ is the dual Hilbert space of $H^1$
 realized through this
embedding. Assume $L^{-1}$ is continuous on $L^{r+1}(\bf{m})$.

For fixed $r>1$, we consider the following Gelfand triple
$$V:=L^{r+1}({\bf m})\subseteq H \subseteq
 V^*$$
and the stochastic porous media equation
 \beq\label{porous media}
  \d X_t^\vv = (L\Psi(t,X_t^\vv)+\Phi(t,X_t^\vv))\d t+\vv
 B(t,X_t^\vv)\d
  W_t, \ X_0^\vv=x\in H.
  \end{equation}
where $W_t$ is a cylindrical
 Wiener process on $L^2(\bf m)$,
$\Psi,\Phi: [0,T]\times \R \to \R$ are measurable and continuous in the second
variable. Suppose $L^2(\mathbf{m})\subseteq H$ is compact and
$B:[0,T]\times V \to L_2(L^2{(\bf m)})$. If  there exist two constants $\delta>0$ and
$K$ such that
\beq\label{condition for porous media equation} \beg{split}
 & |\Psi(t,x)|+|\Phi(t,x)|+\|B(t,0)\|_2 \le
K(1+|x|^{r}),
\ \ t\in [0,T], x\in \R;\\
& -\<\Psi(t,u)-\Psi(t,v),u-v\> -
 \<\Phi(t,u)-\Phi(t,v),L^{-1}(u-v)\>\\
 &\ \ \ \ \ \ \ \le
-\delta\|u-v\|_{V}^{r+1}+K\|u-v\|_H^2; \\
 & \|B(t,u)-B(t,v)\|_{2}^2\le
 K\|u-v\|_H^2,
 \ \ t\in [0,T], u,v\in V.
\end{split}\end{equation}
Then $\{X^\vv\}$ satisfy the large deviation principle on
$C([0,T];H)\cap L^{r+1}([0,T];V)$.

\end{exa}

\beg{proof} From the assumptions and the relation
$$_{V^*}\<L\Phi(t,u)+\Phi(t,u),u\>_V=-\<\Phi(t,u),u\>-\<\Phi(t,u),L^{-1}u\>,$$
it's easy to show that $(A1)-(A4)$
 hold for $\alpha=r+1$ from (\ref{condition for porous media equation}).
  We refer to \cite[Example 4.1.11]{R} for details,
  see also \cite{DRRW,RWW,Wang}.
Hence the conclusion follows from Theorem \ref{T1.1}.
 \end{proof}

\beg{rem} (i) If we take $L$ the Laplace operator on a smooth
bounded domain in a complete Riemannian manifold with Dirichlet
boundary condition. A simple example for $\Psi$ and $\Phi$ satisfy
$(\ref{condition for porous media equation})$ is given by
$$\Psi(t,x)=f(t)|x|^{r-1}x, \ \
\Phi(t,x)=g(t)x$$ for some strictly positive continuous function $f$
and bounded function $g$ on $[0,T]$.

(ii) This example generalized the main result in \cite[Theorem
1.1]{RWW} where $LDP$ was obtained for stochastic porous media equations
with additive noise. In \cite{RWW} the authors mainly used the
piecewise linear approximation to the path of Wiener process and
generalized contraction principle.\end{rem}

If we assume $0<r<1$ in the above example (cf.\cite{LW,RRW}), then the equation is the stochastic version of classical fast diffusion equation.
The behavior of the solutions to these two types of PDE has many
essentially different aspects, see e.g.\cite{Aronson}.
\beg{exa}\label{E6.3} (Stochastic fast diffusion equations)\\
Assume the same framework as Example $\ref{E6.2}$ for $0<r<1$, i.e.
assume the embedding $V:=L^{r+1}({\bf m})\subseteq H$ is continuous
and dense. We consider the equation \beq\label{fast diffusion}
 \d X_t^\vv = \big\{L\Psi(t,X_t^\vv)+ \eta_t X_t^\vv\big\}\d
t +\vv B(t, X_t^\vv) \d W_t,\ X_0^\vv=x\in H,
\end{equation}
where $\eta: [0,T]\to \R$  is locally bounded and measurable and
$\Psi: [0,T]\times \R\to \R$
is  measurable and continuous in the second variable,   $W_t$ is a cylindrical
 Wiener process on $L^2(\bf m)$ and
$B:[0,T]\times
V \to L_2(L^2({\bf m}))$ are measurable.

Suppose there exist constants $\delta>0$ and $K$
 such that for all $x, y\in\R, t\in[0,T]$ and $u,v\in V$
\beq\label{condition for fast diffusion} \beg{split} &
|\Psi(t,x)|+\|B(t,0)\|_2\le K
(1+|x|^r);\\
& (\Psi(t,x)-\Psi(t,y))(x-y)\ge \delta
|x-y|^2(|x|\lor |y|)^{r-1};\\
& \|B(t,u)-B(t,v)\|_{2}^2 \le
 K\|u-v\|_H^2; \\
 & \|B(t,u)\|_{L(L^2({\bf m}), V^*)} \le K(1+\|u\|_V^r).
\end{split}\end{equation}
Then $\{X^\vv\}$ satisfy the large deviation principle on $C([0,T];H)$.
\end{exa}
\begin{proof} Note that
$$ _{V^*}\<L\Psi(t,u)+\eta_t u, u \>_V=- \<\Psi(t,u), u \>_{L^2}+
\<\eta_t u, u \>_H,$$
then it's easy to show $(A1),(A2^\prime),(A3)-(A4)$ hold for
$\alpha=r+1$ under assumptions (\ref{condition for fast diffusion}).
Then the conclusion follows from Theorem \ref{T1.2}.
\end{proof}

\beg{rem}(i) In particular, if $\eta=0, B=0$ and $\Psi(t,s)=
|s|^{r-1}s$ for some $r\in (0,1)$, then (\ref{fast diffusion})
reduces back to the
 classical
fast-diffusion equations (cf. \cite{Aronson}).

(ii) In the example we assume the embedding $L^{r+1}({\bf m})\subseteq H$ is
continuous and dense only for simplicity, see \cite{LW} and \cite[Remark
4.1.15]{R} for
 some sufficient conditions of this assumption. But in general
 $L^{r+1}({\bf m})$ and $H$ are incomparable, hence one need to
 use the more general framework as in \cite{RRW} involving with Orlicz space.
\end{rem}

 \beg{exa}\label{E6.4} (Stochastic $p$-Laplace equation)\\
 Let $\Lambda$ be an open bounded domain in $\mathbb{R}^d$ with smooth boundary. We consider
the  triple
$$V:=H^{1,p}_0(\Lambda)\subseteq H:=L^2(\Lambda)\subseteq
 (H^{1,p}_0(\Lambda))^*$$
and the stochastic $p$-Laplace equation
 \beq\label{p}
 \d X_t^\vv=\left[\mathbf{div}(|\nabla X_t^\vv|^{p-2}\nabla
 X_t^\vv)-\eta_t|X_t^\vv|^{\tilde{p}-2}X_t^\vv\right]\d t+\vv B(t,X_t^\vv)\d
 W_t , X_0^\vv=x\in H,
  \end{equation}
where $2\leq p<\infty, 1\leq\tilde{p}\leq p$, $\eta$ is positive
continuous function and $W_t$ is a cylindrical Wiener process on
$H$. If
$$B(t,v)=\sum_{i=1}^N b_i(v)B_i(t), $$
 where $b_i(\cdot): V\rightarrow \mathbb{R}$ are Lipschitz functions and
 $B_i(\cdot):[0,T]\rightarrow L_2(H)$ are continuous,
 then $\{X^\vv\}$ satisfy the large deviation principle on
$C([0,T];H)\cap L^{p}([0,T];V)$.

 \end{exa}

 \begin{proof} The assumptions for existence and uniqueness of the
 solution were verified in \cite[Example 4.1.9]{R} for $\alpha=p$. Hence we only need
 to prove  $(A2)$ holds here.
 By using (\ref{6.1}) in Lemma \ref{L6.1} we have
 \ce
 && _{V^*}\langle \mathbf{div}(|\nabla u|^{p-2}\nabla
 u)-\mathbf{div}(|\nabla v|^{p-2}\nabla v), u-v
 \rangle_{V}\\
 &=&-\int_\Lambda \langle |\nabla u(x)|^{p-2}\nabla u(x)-|\nabla
 v(x)|^{p-2}\nabla v(x), \nabla u(x)-\nabla
 v(x)
 \rangle_{\R^d} \d x\\
 &\leq&-2^{p-2}\int_\Lambda |\nabla u(x)-\nabla v(x)|^{p}\d x\\
 &\leq& -c\|u-v\|_V^p.
\de where c is a positive constant and follows from the Poincar\'{e}
inequality.

By the monotonicity of function $|x|^{\tilde{p}-2}x$ we know
$$_{V^*}\langle |u|^{\tilde{p}-2}u -|v|^{\tilde{p}-2}v, u-v
 \rangle_{V}\geq 0.$$
 Hence $(A2)$ holds.
 Then the conclusion follows from Theorem
 \ref{T1.1}.
\end{proof}
\beg{rem} If $1<p<2$ in $(\ref{p})$, then the assumption $(A2)$ does not hold. Hence like the case of stochastic
fast diffusion equations, we should apply  Theorem \ref{T1.2} to derive the LDP for $(\ref{p})$ on
$C([0,T]; H)$.
\end{rem}

The following SPDE was studied in \cite{KR,L08a}. The main part of
drift is a high order generalization of the Laplace operator.

 \beg{exa}
 Let $\Lambda$ is an open bounded domain in $\mathbb{R}^1$  and
 $m\in\mathbb{N_+}$, consider
the  triple
$$V:=H^{m,p}_0(\Lambda)\subseteq H:= L^2(\Lambda)\subseteq
 (H^{m,p}_0(\Lambda))^*$$
and the stochastic evolution equation
 \beq \label{sde}\beg{split}
 \d X_t^\vv(x)=&\left[(-1)^{m+1}\frac{\partial}{\partial
 x^m}\left(\left|\frac{\partial^m}{\partial x^m}
 X_t^\vv(x)\right|^{p-2}\frac{\partial^m}{\partial x^m}
 X_t^\vv(x)\right)
 +F(t,X_t(x))\right]\d t\\
 & +\vv B(t,X_t^\vv(x)) \d W_t,  \ \ \  X_0^\vv=x\in H,
 \end{split}
  \end{equation}
where $2\leq p<\infty$, $W_t$ is a cylindrical
 Wiener process on
 $H$  and
 $$F:[0,T]\times
V\rightarrow V^*,\ \  B: [0,T]\times
V\rightarrow L_2(H)$$
are measurable. Suppose $B(t,v)=QB_0(t,v), Q\in L_2(H)$ and
$$\aligned
2{ }_{V^*}\langle F(t,u)-F(t,v), u-v\rangle_{V}& \le
 C\|u-v\|_H^2,\\
 \|B_0(t,u)-B_0(t,v)\|_{L(H)}& \le
 C\|u-v\|_H,\\
\|F(t,u)\|_{V^*}+\|B_0(t,0)\|_{L(H)}& \le
 C(1+\|u\|_V^{p-1}), \ u,v\in V, \ t\in[0,T].
\endaligned$$
where $C$ is a constant.  Then $\{X^\vv\}$ satisfy the large deviation
principle on $C([0,T];H)\cap L^p([0,T];V)$.
 \end{exa}
\begin{proof}
By using Lemma \ref{L6.1}, $(A2)$ can be verified  by the
 same argument as in Example \ref{E6.4}. And $(A1),(A3),(A4)$ follow from the assumptions obviously,
hence the conclusion follows from
 Theorem \ref{T1.1}.

\end{proof}

\section*{Acknowledgements} The author would like to thank Professor
 Michael R\"{o}ckner, Fengyu Wang and Xicheng Zhang for their  valuable discussions, and also  thank
Professor Paul Dupuis and Amarjit Budhiraja for the stimulating
communications on their works. Many helpful comments from the  referee are also gratefully acknowledged.


\begin{thebibliography}{9}

\bibitem{Aronson} D.G. Aronson, \emph{The porous medium equation,} Lecture Notes
in Mathematics 1224, Springer, Berlin, 1--46, 1986.


\bibitem{Az} R.G. Azencott, \emph{Grandes  deviations et applications,
 Ecole d'Et\'e de Probabilit\'es de Saint-Flour VII,}
  Lecture Notes in Mathematics  774, 1980.


\bibitem{Be} A. Bensoussan, \emph{Filtrage optimale des systemes
 lin\'eaires,}
  Dunod, Paris, 1971.


\bibitem{BT} A. Bensoussan and R. Temam, \emph{Equations aux derives
 partielles stochastiques non lin\'eaires,}
  Isr. J. Math. 11(1972), 95--129.





\bibitem{BD} A. Budhiraja and P. Dupuis, \emph{A variational
 representation for positive functionals
of infinite dimensional Brownian motion,}  Probab. Math. Statist.
 20(2000), 39--61.


\bibitem{BDM} A. Budhiraja, P. Dupuis and V. Maroulas, \emph{Large
 deviations for infinite
dimensional stochastic dynamical systems,}  Ann. Probab. 36(2008), 1390-1420.

\bibitem{Br} W. Bryc, \emph{Large deviations by the asymptotic value
 method,}
In M. Pinsky, editor, Diffusion Processes and Related Problems in
 Analysis, vol.1, 447--472.
Birkh\"auser, Boston, 1990.

\bibitem{CR} S. Cerrai and M. R\"ockner, \emph{Large deviations
for stochastic reaction-diffusion systems with multiplicative noise
and non-Lipschitz reaction term,} Ann.  Probab. 32(2004),
1100--1139.

\bibitem{C} P.L. Chow, \emph{Large deviation problem for some parabolic
 It\^{
o} equations,} Commun. Pure Appl. Math. 45(1992), 97--120.


\bibitem{DR1} G. Da Prato and M. R\"ockner, \emph{Weak solutions to
 stochastic
porous media equations,} J. Evolution Equ. 4(2004), 249--271.


\bibitem{DRRW} G. Da Prato, M. R\"ockner, Rozovskii and F.-Y.
Wang, \emph{Strong solutions to stochastic generalized porous media
equations: existence, uniqueness and ergodicity,}  Comm. Part. Diff.
Equat.  31(2006), no.2, 277-291.

\bibitem{DaZa} G. Da Prato and J. Zabczyk, \emph{ Stochastic Equations
 in
Infinite Dimensions,} Encyclopedia of Mathematics and its
Applications, Cambridge University Press. 1992.

\bibitem{DZ} A. Dembo and O. Zeitouni,  \emph{Large deviations
techniques and applications,} Springer-Verlag, New
 York. 2000.

\bibitem{DV} M.D. Donsker and S.R.S. Varadhan, \emph{Asymptotic
 evalution of certain Markov
process expectations for large time, I, II, III,} Comm. Pure Appl.
 Math. 28(1975), 1--47;
28(1975), 279--301; 29(1977), 389--461.

\bibitem{DM08} J. Duan and A. Millet,
\emph{Large deviations for the Boussinesq Equations under Random Influences,}
Preprint.

\bibitem{DE} P. Dupuis and R. Ellis, \emph{A weak convergence
 approach
 to the theory of
large deviations,} Wiley, New York. 1997.




\bibitem{FK} J. Feng and T.G. Kurtz, \emph{Large Deviations of Stochastic Processes,}
vol. 131 of \emph{Mathematical Surveys and Monographs,} American
Mathematical Society, Providence, RI, 2006.


\bibitem{Fr} M.I. Freidlin, \emph{Random perturbations of
 reaction-diffusion equations: the quasi-deterministic approximations,}
 Trans. Amer.
 Math. Soc. 305 (1988), 665--697.


\bibitem{FW} M.I. Freidlin and A.D. Wentzell, \emph{Random
 perturbations of dynamical systems,} Translated from the Russian by
 Joseph Szu"cs.
 Grundlehren der Mathematischen Wissenschaften [Fundamental Principles
 of Mathematical Sciences], 260. Springer-Verlag, New York, 1984.

\bibitem{GM} I. Gy\"ongy and A. Millet, \emph{On discretization schemes
 for
stochastic evolution equations,} Pot. Anal.  23(2005), 99--134.



\bibitem{KR} N.V. Krylov and B.L. Rozovskii, \emph{Stochastic evolution
 equations,}
Translated from Itogi Naukii Tekhniki, Seriya Sovremennye Problemy
Matematiki 14(1979), 71--146, Plenum Publishing Corp. 1981.

\bibitem{L08a} W. Liu, \emph{Harnack inequality and applications
 for stochastic evolution equations  with monotone drifts,}
 SFB-Preprint 09-023.

\bibitem{LW} W. Liu and F.-Y. Wang, \emph{Harnack inequality and strong
 Feller
property for stochastic fast diffusion equations,}
  J. Math. Anal. Appl. 342(2008), 651-662.

\bibitem{Pa} E. Pardoux, \emph{Equations aux d\'eriv\'ees partielles
 stochastiques non lin\'eaires monotones,}
Thesis, Universit\'e  Paris XI, 1975.


\bibitem{P} S. Peszat, \emph{Large deviation principle for stochastic
evolution equations,} Probab. Theory Relat. Fields. 98(1994),
113--136.

\bibitem{Pu} A.A. Pukhalskii,
\emph{On the theory of large deviations,} Theory probab. Appl.
 38(1993), 490--497.



\bibitem{RRW} J. Ren, M. R\"ockner and F.-Y. Wang,
\emph{Stochastic generalized porous media and fast diffusion
equations,} J. Diff. Equat. 238(2007), 118--152.

\bibitem{RZ} J. Ren and X. Zhang, \emph{Freidlin-Wentzell large
 deviations for homeomorphism flows of
non-Lipschitz SDE,} Bull. Sci.  129(2005), 643--655.

\bibitem{RZ1} J. Ren and X. Zhang, \emph{Schilder theorem for the
 Brownian motion on
 the diffeomorphism group of the circle,} J. Funct. Anal.  224(2005),
107--133.

\bibitem{R} C. Pr\'{e}v\^{o}t and M. R\"ockner, \emph{A Concise Course
 on Stochastic Partial Differential Equations,} Lecture Notes in
Mathematics  1905, Springer, 2007.



\bibitem{RSZ}  M. R\"{o}ckner, B. Schmuland and X. Zhang,
\emph{Yamada-Watanabe Theorem for stochastic evolution equations in
infinite dimensions,} Condensed Matter Physics  54(2008), 247-259.

\bibitem{RWW} M. R\"ockner, F.-Y. Wang and L. Wu, \emph{Large
 deviations
for stochastic Generalized Porous Media Equations ,} Stoch. Proc.
Appl.
 116(2006), 1677--1689.

\bibitem{SS} S.S. Sritharan and P. Sundar, \emph{Large deviations for
 the two-dimensional Navier-Stokes
 equations with multiplicative noise,} Stoc. Proc. Appl.  116(2006),
1636--1659.

\bibitem{St} D.W. Stroock, \emph{An Introduction to the Theory of Large
 Deviations,} Spring-Verlag, New York, 1984.

\bibitem{Va1} S.R.S. Varadhan, \emph{Asymptotic probabilities and
 differential equations,}
  Comm. Pure Appl. Math.  19 (1966), 261--286.


\bibitem{Va2} S.R.S. Varadhan, \emph{Diffusion processes in a small
 time interval,}
  Comm. Pure Appl. Math.  20  (1967), 659--685.


\bibitem{Va} S.R.S. Varadhan, \emph{Large Deviations and Applications,}
 CBMS 46, SIAM, Philadelphia, 1984.


\bibitem{Wa} J.B. Walsh, \emph{An introduction to stochastic partial
 differential equations,}
  Ecole d'Ete de Probabilite de Saint-Flour XIV (1984), P.L. Hennequin
 editor, Lecture Notes in Mathematics 1180
,  265--439.

\bibitem{Wang} F.-Y. Wang, \emph{Harnack Inequality and Applications for Stochastic Generalized Porous Media equations,}
 Ann. Probab.  35(2007), 1333-1350.


\bibitem{Wu04} L. Wu,
\emph{On large deviations for moving average processes, } In
Probability, Finance and Insurance, pp.15-49, the proceeding of a
Workshop at the University of Hong-Kong (15-17 July 2002), Eds: T.L.
Lai, H.L. Yang and S.P. Yung. World Scientific 2004, Singapour.

\bibitem{Z} E. Zeidler, \emph{Nonlinear Functional Analysis and
its Applications, II/B, Nonlinear Monotone Operators,}
Springer-Verlag, New York: 1990.

\bibitem{Zh} X. Zhang, \emph{On Stochastic evolution equations with
non-Lipschitz coefficients,} Preprint.

\end{thebibliography}
\end{document}